\documentclass{article}

\usepackage{arxiv}

\usepackage[utf8]{inputenc}
\usepackage[T1]{fontenc}
\usepackage{microtype}

\usepackage{amsfonts}
\usepackage{amsmath}
\usepackage{amssymb}
\usepackage{amsopn}
\usepackage{amsthm}
\usepackage{nicefrac}

\usepackage{booktabs}
\usepackage{graphicx}
\usepackage{epstopdf}
\usepackage{makecell}
\usepackage{enumitem}
\usepackage{subcaption}

\usepackage{algorithm}
\usepackage{algpseudocode}

\theoremstyle{plain}

\newtheorem{theorem}{Theorem}[section]
\newtheorem{lemma}{Lemma}[section]
\newtheorem{proposition}{Proposition}[section]

\theoremstyle{definition}

\newtheorem{definition}{Definition}[section]

\theoremstyle{remark}

\newtheorem{remark}{Remark}[section]

\usepackage{cite}
\usepackage{url}
\usepackage{hyperref}
\usepackage{cleveref}

\crefname{theorem}{theorem}{theorems}
\Crefname{theorem}{Theorem}{Theorems}

\crefname{lemma}{lemma}{lemmas}
\Crefname{lemma}{Lemma}{Lemmas}

\crefname{proposition}{proposition}{propositions}
\Crefname{proposition}{Proposition}{Propositions}

\crefname{definition}{definition}{definitions}
\Crefname{definition}{Definition}{Definitions}

\crefname{assumption}{assumption}{assumptions}
\Crefname{assumption}{Assumption}{Assumptions}

\crefname{remark}{remark}{remarks}
\Crefname{remark}{Remark}{Remarks}

\crefname{example}{example}{examples}
\Crefname{example}{Example}{Examples}

\crefname{hypothesis}{hypothesis}{hypotheses}
\Crefname{hypothesis}{Hypothesis}{Hypotheses}

\crefname{fact}{fact}{facts}
\Crefname{fact}{Fact}{Facts}

\setlist[enumerate]{leftmargin=.5in}
\setlist[itemize]{leftmargin=.5in}

\renewcommand{\headeright}{}
\renewcommand{\undertitle}{}
\renewcommand{\shorttitle}{Learning piecewise-smooth dynamical systems}

\newcommand{\p}{\mathcal{P}}
\newcommand{\tp}{\widetilde{\p}}
\newcommand{\R}{\mathbb{R}}

\title{Learning piecewise-smooth dynamical systems}

\author{
  \textbf{Davide Murari} \\
  Department of Applied Mathematics\\
  and Theoretical Physics \\
  University of Cambridge \\
  \texttt{dm2011@cam.ac.uk} \\
  \And
  \textbf{Erik Jansson} \\
  Department of Applied Mathematics\\
  and Theoretical Physics \\
  University of Cambridge \\
  \texttt{eoj23@cam.ac.uk} \\
  \AND
  \textbf{Chris Budd OBE} \\
  Department of Mathematical Sciences \\
  University of Bath \\
  \texttt{mascjb@bath.ac.uk} \\
  \And
  \textbf{Carola-Bibiane Schönlieb} \\
  Department of Applied Mathematics\\
  and Theoretical Physics \\
  University of Cambridge \\
  \texttt{cbs31@cam.ac.uk} \\
}

\date{}

\hypersetup{
    hidelinks,
    pdftitle={Learning piecewise-smooth dynamical systems},
    pdfauthor={Davide Murari, Erik Jansson, Chris Budd OBE, Carola-Bibiane Schönlieb}
}

\begin{document}

\maketitle

\begin{abstract}
Discovering dynamical systems from trajectory data is a central problem in applied mathematics and engineering. Whilst recent advances in machine learning have led to strong progress in data-driven system identification, much less attention has been given to systems with discontinuous dynamics. These systems are nevertheless highly relevant in applications, including climate dynamics and mechanical systems with friction. In this work, we consider the problem of identifying piecewise-smooth dynamical systems directly from trajectory data. Compared with the smooth setting, this requires recovering the governing equations and detecting the switching hyperplanes that separate different dynamical regimes and characterising their behaviour, such as sliding motion. We present a modular framework for discovering such systems by first estimating switching hyperplanes from data and then learning smooth dynamics within each region using geometry-constrained neural networks. The geometry-learning phase is studied from a statistical perspective, analysing the identifiability of the discontinuities and the robustness of the procedure. We also introduce a novel neural network architecture with a prescribed discontinuity set, and provide a theoretical analysis of its approximation properties. The approach is tested on low-dimensional benchmark problems, including dry-friction oscillators and the PP04 climate model for the ice ages.
\end{abstract}


\section{Introduction}\label{sec:intro}
Push a heavy box across a floor, and it does not move at once, until suddenly, it does.
Open a stiff door slowly and the hinge does not turn smoothly but in tiny catches and releases.
In both cases, the dynamics do not vary smoothly but jump between distinct regimes, sticking and slipping, as the system passes from one regime to another.

Systems like these are not unusual. We consider \emph{piecewise-smooth systems}, whose vector fields are smooth within separate regions $\mathcal P_k$ of state space but may jump across the \emph{switching hyperplanes} $\Sigma_{k,\ell}$ dividing them.
Throughout this work, such systems are interpreted in the sense of Filippov, which replaces the discontinuous vector field on a switching hyperplane by a set-valued convexification of its nearby limiting values. We refer to systems equipped with this solution concept as \emph{Filippov systems}. At a regular switching hyperplane, Filippov solutions may cross the hyperplane or slide along an attracting part of it \cite{filippov1988differential,Laurea2008book,jeffreyModelingBook,brogliato2016nonsmooth}. Discontinuous differential equations arise in mechanical systems with friction \cite{Karnopp1985}, sliding-mode control \cite{Edwards1998}, robotics \cite{Arkin1998-vk}, and the glacial cycles of the Earth's climate \cite{Paillard2004, Budd2022, Morupisi2020}.

\begin{figure}[ht!]
    \centering
    \includegraphics[width=0.8\linewidth]{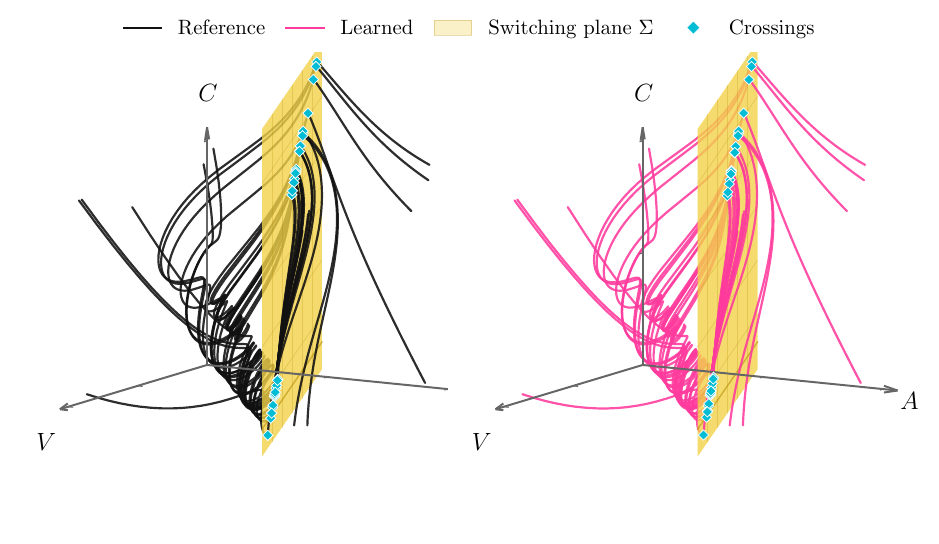}
    \caption{The PP04 system \cite{Budd2022,Morupisi2020} is a
    three-dimensional climate model of the ice ages. The vector field is
    discontinuous across the affine plane $\Sigma$, shown in yellow, and the cyan markers indicate trajectory crossings of this plane. The left panel shows trajectories of the reference system, while the right panel shows
    rollouts of the learned discontinuous model from the same initial conditions. The learned model is supplied with the true switching plane. The close agreement between the two trajectory families demonstrates its ability to reproduce the discontinuous dynamics. Further details are provided in \Cref{sec:traj-pp04_experiment}.}
    \label{fig:pp04-orbits}
\end{figure}

Often the governing dynamics are unknown and must be recovered from observed trajectories. For piecewise-smooth systems, this requires learning the vector field in each regime, locating the switching hyperplanes, and characterising the behaviour at those hyperplanes. A misplaced switching hyperplane misplaces the discontinuities, regardless of how accurately the local dynamics are learned. \Cref{fig:pp04-orbits} illustrates this setting for the PP04 climate model: the discontinuity is confined to an affine plane, and the learned model reproduces reference trajectories that repeatedly cross this plane.

Data-driven modelling of dynamical systems is commonly pursued in two complementary ways.  
One seeks to recover explicit governing equations from observations, deriving a symbolic approximation of them \cite{brunton2016discovering,cranmer2020discovering,cranmer2023pysr}. 
A related difficulty is that simple governing equations may only exist in appropriate coordinates. This has motivated autoencoder-based approaches that simultaneously learn latent coordinates and parsimonious dynamics \cite{lusch2018deep,champion2019data}. 
Another way to model dynamical systems instead learns the time evolution directly from trajectory data using neural or latent dynamical models.
This includes recurrent architectures for time-series forecasting \cite{hochreiter1997long,vlachas2018data}, neural ordinary differential equations \cite{chen2018neural}, and latent or controlled variants designed for partially observed or irregularly sampled trajectories \cite{rubanova2019latent,kidger2020neuralcde}. 
Universal differential equations further combine mechanistic model components with neural parameterisations of unknown terms \cite{rackauckas2020universal}. 

Within scientific machine learning, an important refinement is to constrain the architecture so that the learned model respects the known structure of the dynamics. 
For Hamiltonian and Lagrangian systems, this has led to Hamiltonian neural networks, Lagrangian neural networks, symplectic networks, and symplectic neural flows, where conservation laws, constraints, or symplecticity are built into the parametrisation \cite{greydanus2019hamiltonian,cranmer2020lagrangian,jin2020sympnets,celledoni2023learning,canizares2024symplectic}. 
The present work follows the same structure-preserving philosophy but applies it to non-smooth systems. The relevant structure in our setting is the geometry of the switching set and the requirement that learned discontinuities remain confined to it. Existing work has considered learning local dynamics when the switching geometry is known, including friction problems \cite{lathourakis2024physics}. Here, we additionally seek to recover that geometry from trajectory data.

\paragraph{Main contributions}
We introduce a geometry-first, two-stage framework for discovering
piecewise-smooth dynamical systems from trajectory data. Our main contributions are three.

\begin{figure}[ht!]
    \centering
    \includegraphics[width=0.6\linewidth]{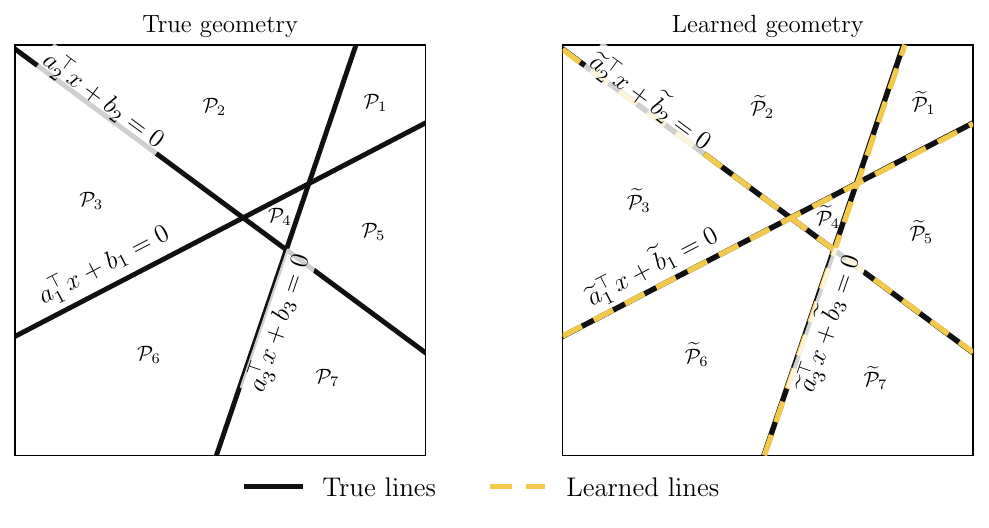}
    \caption{Schematic illustration of the geometry-learning stage.  The true
    switching geometry is an arrangement of affine lines, which partitions the
    state space into polyhedral regions $\mathcal P_i$.  From trajectory data, the
    algorithm estimates the switching geometry and obtains an approximate arrangement,
    shown on the right, with corresponding regions $\widetilde{\mathcal P}_i$. More details on this experiment are in Appendix \ref{app:three-line-appendix}.}
    \label{fig:poly}
\end{figure}

First, we develop a method for recovering arrangements of affine switching hyperplanes before fitting the vector field. Their sign patterns define interpretable polyhedral regions, see \Cref{fig:poly}, so the inferred switching geometry can be inspected and validated independently of the learned dynamics. We complement this method with a detectability analysis quantifying how jump size, observation noise, and derivative-estimation error affect geometry recovery.

Second, we introduce a geometry-preserving neural architecture whose discontinuities are confined, by construction, to a prescribed hyperplane arrangement. It can therefore learn distinct dynamics in adjacent regions without introducing discontinuities elsewhere. The regression experiments in \Cref{sec:regression} demonstrate this distinction directly.

Third, we prove that this architecture can approximate piecewise-smooth vector fields to arbitrary accuracy on compact sets away from the switching interfaces, while preserving the prescribed discontinuity set. Experiments on dry friction oscillators and the PP04 climate model validate both stages of the framework. In the reported baseline comparisons, the proposed models achieve lower rollout and vector-field errors than smooth Neural ODEs while using an order of magnitude fewer trainable parameters.

\paragraph{Modularity and interpretability of the methodology}
The two stages of the framework are modular. The geometry-estimation stage is driven by {\em derivative approximations} from trajectory data, and the derivative-estimation procedure can be adapted to the noise level, sampling rate, and dimension of the problem.  In this paper, this stage combines derivative-jump candidate extraction with random sample consensus (RANSAC) hyperplane fitting, followed by total least squares refinement. This stage can be skipped when the switching geometry is known in advance. 
The dynamics-learning stage can likewise incorporate prior knowledge about the local dynamics. 
For example, if each smooth piece is known to be affine, one can use a piecewise-affine model. If the local dynamics are Hamiltonian, one can use a Hamiltonian neural network \cite{greydanus2019hamiltonian} in each region. 
In our numerical experiments, we do not assume such detailed prior knowledge. Instead, we use a general piecewise-smooth neural vector field. 

\paragraph{Outline of the paper} This paper is structured as follows: The necessary background on piecewise-smooth dynamical systems is presented in \Cref{sec:background}. \Cref{sec:learning} describes the dataset and the joint learning of geometry and dynamics. Some theoretical aspects of the proposed methodology are considered in \Cref{sec:theory}, including an analysis of the hyperplane detectability from noisy data and the approximation properties of the considered network. We include a selection of numerical experiments in \Cref{sec:experiments} looking at friction and climate models. These experiments illustrate both the methodology and the theoretical results. Finally, we draw some conclusions from this work. 

\section{Piecewise-smooth dynamical systems}
\label{sec:background}
In this work, we use the term \textit{piecewise-smooth system} to refer to a dynamical system whose vector field may be discontinuous across affine switching hyperplanes $\Sigma_{k,\ell}$ \cite{filippov1988differential} and is continuously differentiable away from them. Its solutions are understood in the Filippov sense \cite{filippov1988differential,Laurea2008book}. We now make this setting precise. 

We focus on systems with smoothness regions given by polyhedral partitions induced by finite arrangements of affine hyperplanes. More precisely, let $\{H_1,\ldots,H_R\}$ be a finite affine hyperplane arrangement of $\R^d$, and let
\(P=\{\p_1,\ldots,\p_K\}\) be the collection of connected components of 
\[
\mathbb R^d\setminus \bigcup_{r=1}^R H_r.
\]
Then each $\p_k$, $k=1,\dots,K$, is an open convex cell of the form
\begin{equation}\label{eq:partitioning}
\p_k=\{x\in\R^d\mid A_k x<b_k\},
\qquad
A_k\in\R^{m_k\times d},\; b_k\in\R^{m_k},
\end{equation}
where the inequality is understood componentwise, and such that $\bigcup_{k=1}^K \overline{\p_k} = \R^d$, where $\overline{\p_k}$ denotes the closure of $\p_k$. For every pair of regions sharing a codimension-one boundary, we denote by $\Sigma_{k,\ell}:=\operatorname{aff}(\overline{\p_k}\cap\overline{\p_\ell})$ the affine hyperplane containing the facet shared by $\p_k$ and $\p_\ell$. We refer to $\Sigma_{k,\ell}$ as a \emph{switching hyperplane}.
We consider dynamical systems governed by piecewise-$\mathcal C^1$ vector fields $f: \R^d\times\R\to\R^d$, namely systems of the form
\begin{equation}\label{eq:main-system_na}
    \dot x = f(x,t) := f_k(x,t), \quad x \in \p_k, \quad k = 1,\dots,K,
\end{equation}
where each $f_k : \p_k\times \R\to \R^d$ is continuously differentiable. Thus, all discontinuities of $f$ are confined to the shared facets, each of which lies in its corresponding switching hyperplane $\Sigma_{k,\ell}$.
We denote the union of these facets by $\Sigma$, that is,
\[
\Sigma := \R^d\setminus \bigcup_{k=1}^K \p_k = \bigcup_{r=1}^R H_r.
\]
Throughout this work, we assume that each $f_k$ admits a continuously differentiable extension $f_k:\R^d\times\R\to\R^d$, and we use the same notation for the vector field on $\p_k$ and its extension.

The values of $f$ on $\Sigma$ may be assigned arbitrarily without changing its Filippov regularisation: because the facets have Lebesgue measure zero, the regularisation depends on the neighbouring limiting vector fields rather than on pointwise values assigned on the facets. At regular points of codimension-one facets, this regularisation simplifies, as discussed in more detail below. At intersections of several facets, the regularisation accounts for the limiting vector fields from all adjacent regions. A Filippov solution is an absolutely continuous curve whose derivative belongs to the Filippov regularisation of $f$ almost everywhere \cite{filippov1988differential,Laurea2008book}.

A typical example is the case of two half-spaces separated by a hyperplane
\[
  \Sigma:=\left\{x\in\mathbb{R}^d:\ s(x)=0\right\},
  \qquad s(x):=a^\top x+b,
\]
where $a\in\mathbb{R}^d$ satisfies $\lVert a\rVert_2=1$ and $b\in\mathbb{R}$, so that
\begin{equation}\label{eq:main-system}
  \dot x = f(x,t):=
  \begin{cases}
    f^+(x,t), & s(x)>0,\\[1mm]
    f^-(x,t), & s(x)<0.
  \end{cases}
\end{equation}
This single-hyperplane setting describes the local behaviour near a regular codimension-one facet separating two regions of a more general partition. For $x\in\Sigma$, the normal components of the two limiting vector fields are
\[
\nu^+(x,t):=a^\top f^+(x,t),
\qquad
\nu^-(x,t):=a^\top f^-(x,t).
\]
If $\nu^+$ and $\nu^-$ have the same sign, trajectories locally cross the switching hyperplane. If
\[
\nu^+(x,t)<0<\nu^-(x,t),
\]
then both limiting vector fields point towards $\Sigma$, and the hyperplane is locally attracting at $(x,t)$. In the Filippov convention \cite{filippov1988differential}, {\em sliding motion} on the locally attracting part of the hyperplane $\Sigma$ is then described by the convex combination
\[
f_{\Sigma}(x,t)
=
\alpha(x,t)f^+(x,t)
+
(1-\alpha(x,t))f^-(x,t),
\]
where $\alpha\in[0,1]$ is chosen so that the trajectory remains tangent to the hyperplane, i.e., 
\[
a^\top f_{\Sigma}(x,t)=0.
\]
Equivalently,
\[
\alpha(x,t)
=
\frac{\nu^-(x,t)}{\nu^-(x,t)-\nu^+(x,t)}.
\]
This gives the sliding vector field along $\Sigma$. The opposite sign
configuration,
\[
\nu^-(x,t)<0<\nu^+(x,t),
\]
corresponds to a locally repelling or escaping part of the hyperplane, where forward solutions starting on $\Sigma$ may fail to be unique.
\begin{figure}
    \centering
    \includegraphics[width=.7\linewidth]{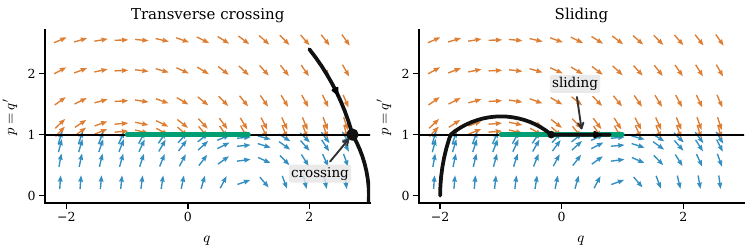}
    \caption{Basic switching behaviours in the toy Filippov system $\dot q=p,\ \dot p=-q+\operatorname{sign}(1-p)$. The switching hyperplane is $\Sigma=\{p=1\}$, shown by the horizontal black line. The green segment marks the attracting sliding region $-1<q<1$.  On this segment, the one-sided normal velocities satisfy $-q-1<0<-q+1$, so the Filippov vector field is tangent to $\Sigma$ and gives $\dot q=1,\ \dot p=0$.  The black curves show representative trajectories: one crosses $\Sigma$ transversally outside the sliding interval, while the other reaches the attracting part of $\Sigma$ and then slides along it.
}
\label{fig:switching-geometry}
\end{figure}
\Cref{fig:switching-geometry} illustrates the two behaviours that are most important for the learning problem: transversal crossings, which generate localised derivative jumps, and attracting sliding, where trajectories remain on the switching hyperplane and must be propagated with the Filippov sliding vector field $f_\Sigma$.

For further background, see \cite{Laurea2008book,filippov1988differential} and the references therein.

\section{Learning piecewise-smooth dynamical systems}\label{sec:learning}

Our approach to learning piecewise-smooth dynamical systems is split into two sequential tasks.

\noindent {\em First}, we estimate the switching geometry by recovering an arrangement of affine
switching hyperplanes. The sign patterns of this arrangement define the
approximate polyhedral regions
\[
\widetilde P=\{\tp_1,\ldots,\tp_{\widetilde{K}}\},
\]
which are used during dynamics learning.

\noindent {\em Second}, we learn a piecewise-smooth vector field adapted to this recovered partition. 

\noindent This procedure separates ``geometry learning'' from ``dynamics learning'', and allows inductive biases specific to each task to be introduced independently.

We now describe the core building blocks of our methodology, starting with the dataset, then moving to the geometry and dynamics learning phases.

\subsection{Dataset}\label{sec:dataset}

We consider $N$ reference trajectories and their noisy observations, indexed by $n=1,\ldots,N$. The $n$th trajectory starts from the initial condition $x_0^n$ at time $t_0^n$, and is sampled at the regularly spaced times
\[
t_j^n:=t_0^n+j \; h,\qquad j=0,\ldots,J,
\]
with final time $T_n=t_J^n=t_0^n+Jh$. For simplicity, all trajectory segments are taken to have the same number of time steps. The extension to segments of non-uniform length is immediate. Let $x(t;t_0,x_0)$ denote the solution satisfying
\[
\dot x(t;t_0,x_0)=f(x(t;t_0,x_0),t),
\qquad
x(t_0;t_0,x_0)=x_0.
\]
For the $n$th trajectory, we write
\[
x^n_*(t)=x(t;t_0^n,x_0^n),
\qquad t\in[t_0^n,T_n].
\]
We assume that the sampling times are known exactly and are not subject to timing errors. In the idealised data model, the noise-free states $x_*^n(t_j^n)$ are additively perturbed by independent Gaussian observation noise. The observations are thus
\begin{equation}
\label{eq:data}
  y_j^n:=x^n_*(t_j^n)+\eta_j^n,
  \qquad j=0,\ldots,J,
\end{equation}
where $\eta_j^n\sim\mathsf{N}(0,\sigma^2 I_d)$ are independent across both $n$ and $j$, $I_d\in\R^{d\times d}$ denotes the identity matrix, and $\sigma\ge 0$ is the noise level. For numerically generated data, the discrepancy from the exact solution may also contain time-discretisation errors. The analysis in \Cref{sec:theory_geometryLearning} isolates the effect of Gaussian observation noise. In the numerical experiments, see \Cref{sec:experiments}, the reference
trajectories are generated synthetically from the benchmark systems. When an exact solution is available between switching events, we use it together with event localisation at switching times; otherwise, we use an event-aware numerical solver. The system-specific sampling boxes, time grids, and solver details are given in \Cref{sec:trajectory_data_generation}.

\subsection{Learning the geometry}\label{sec:geom_learning}

There are two natural signals for recovering $\Sigma$: (i) trajectory points lying on sliding segments, and (ii) jumps in estimated velocities across consecutive samples. The first signal is present only when sliding occurs and is highly sensitive to noise, which can break collinearity. We therefore base the geometry-learning stage on the second signal: large changes in estimated velocity are treated as candidate discontinuity points.

Implementing velocity-jump-based geometry learning requires choosing a velocity estimator. In principle, our method is independent of this choice, but its performance depends on the properties of the adopted approach. 
For instance, in the regime of negligible observational noise, we can compute velocity estimates along trajectory $n$ at the points $y^n_j$ using central finite differences, i.e., 
\begin{align}
\label{eq:central_diff}
v_{h,j}^n = \frac{y^n_{j+1}-y^n_{j-1}}{2h},
\qquad j=1,\ldots,J-1,
\end{align}
meaning that each trajectory provides $J-1$ finite-difference estimates. 

If noise is present, we can instead employ a noise-mitigating Savitzky--Golay velocity estimator \cite{VanBreugel2020, Savitzky1964}, so that  
\begin{align}
    \label{eq:savgol}
    v_{h,j}^n = \sum_{k = -r}^r s_k \; y^n_{j+k},
    \qquad j=r,\ldots,J-r,
\end{align}
where the coefficients are obtained by fitting a polynomial of degree $p$ to uniformly spaced coordinate indices by linear least squares and therefore do not depend on the observed data. Nonlinear alternatives include total-variation regularised differentiation \cite{VanBreugel2020, Chartrand2011}. Our numerical tests use Savitzky--Golay differentiation as a shared geometry-learning procedure across the experiments because it provides accurate hyperplane recovery at modest computational cost.

\noindent We then compute the velocity jump $\Delta v_{h,j}^n = v_{h,j+1}^{n}-v_{h,j}^{n}$. Note that \Cref{eq:central_diff} and \Cref{eq:savgol} are \emph{linear} estimators, that is, they are of the form 
\begin{align}
    \label{eq:lin_est}
     v_{h,j}^n  = \sum_{k=-r}^r w_k \; y^n_{j+k}
\end{align}
for some cut-off $r \in \mathbb{N}$ and weights $w_k\in\mathbb{R}$, $k=-r,\ldots,r$. 
Linear estimators are inexpensive and easy to analyse, but can be sensitive to observational noise, as illustrated in \Cref{sec:learn-geometry}.
\medskip

\noindent\textbf{Algorithm overview.}
\Cref{alg:single_hyperplane_recovery_experiments} states the geometry-learning procedure for a single dominant affine switching hyperplane, which is the setting used in the main experiments. Having
computed the jumps, we rank all the derivative-jump magnitudes globally and retain the largest fraction $\rho\in(0,1)$. For each retained jump interval, we also include neighbouring intervals, since a true crossing may affect several adjacent velocity estimates. The endpoints of these intervals form an overcomplete candidate cloud $\mathcal C$. This cloud contains points near the switching hyperplane, but also outliers caused by curvature, noise, and derivative-estimation artefacts. We proceed with the hyperplane recovery only if there are enough candidates, i.e., $|\mathcal C| \geq \max\{3,m_{\min},d\}$ with $m_{\min}=30$ for the single-hyperplane recovery tasks. We therefore fit the hyperplane using random sample consensus (RANSAC) \cite{Fischler1981}, a robust randomised fitting procedure that repeatedly proposes a hyperplane from a minimal subset of points, counts the number of candidate points within a prescribed distance tolerance, and keeps the hyperplane with the largest consensus set.  We then refine the selected inliers by total least squares.  See Appendix~\ref{app:ransaclinreg} for more background on RANSAC.

\begin{remark}
\noindent This stage requires the observed trajectories to produce velocity jumps that dominate smooth curvature and derivative-estimation error. Transversal crossings and transitions into or out of sliding can provide such a signal, whereas grazing contacts and trajectories confined to one side of the switching hyperplane may not. For example, in $\dot x=\mathrm{sign}(y)$, $\dot y=0$, the vector field is discontinuous across $\Sigma=\{y=0\}$, but every trajectory has constant $y$ and therefore produces no derivative jump. Thus, identifiability requires both a sufficiently large jump and sufficient interaction between the observed trajectories and $\Sigma$. We analyse the jump-to-noise balance in \Cref{sec:theory_geometryLearning} and test it numerically in \Cref{sec:experiments}.
\end{remark}

\begin{algorithm}[ht!]
\caption{Single-hyperplane recovery}
\label{alg:single_hyperplane_recovery_experiments}
\begin{algorithmic}[1]
\Require Warmup trajectories $\{y^n_j\}_{n=1,j=0}^{N,J}\subset\mathbb{R}^d$; step size $h$; Savitzky--Golay window $w$ and polynomial degree $p$; retained jump fraction $\rho\in(0,1)$ and neighbouring interval radius $R$; RANSAC tolerance $\varepsilon$; RANSAC iterations $M$; minimum support $m_{\min}$.
\Ensure Unit normal $a\in\mathbb{R}^d$ and offset $b\in\mathbb{R}$ defining the recovered switching hyperplane $\Gamma=\{x\in\mathbb{R}^d:a^\top x+b=0\}$.
\State Estimate velocities with Savitzky--Golay differentiation:
\Statex\hspace{1.5em}$\widehat v^n_j\gets \operatorname{SG}_{w,p}[y^n]_j$, for $n=1,\dots,N$ and $j=0,\dots,J$.
\State Compute derivative-jump magnitudes
\Statex\hspace{1.5em}$J_n(j)\gets \|\widehat v^n_{j+1}-\widehat v^n_j\|_2$, for $j=0,\dots,J-1$.
\State Rank all jump intervals $(n,j)$ globally by $J_n(j)$.
\State Set $L\gets NJ$.
\State Keep the largest $\max\{1,\operatorname{round}(\rho L)\}$ jump intervals.
\State Add the immediate neighbouring intervals of each retained jump:
\Statex\hspace{1.5em}$\mathcal E\gets\{(n,\ell):\ell\in\{j-R,\dots,j,\dots,j+R\},\ 0\leq \ell<J,\ (n,j)\text{ is retained}\}$.
\State Build the candidate cloud from both endpoints of every selected interval:
\Statex\hspace{1.5em}$\mathcal C\gets\{y^n_\ell,y^n_{\ell+1}:(n,\ell)\in\mathcal E\}$.
\State Remove duplicate points from $\mathcal C$.
\State Set $m_{\rm req}\gets\max\{3,m_{\min},d\}$.
\If{$|\mathcal C|<m_{\rm req}$}
    \State \textbf{fail}.
\EndIf
\State $(a_\star,b_\star,\mathcal I_\star)
\gets
\operatorname{RANSAC}(\mathcal C,\varepsilon,M,m_{\rm req})$.
\Comment{For a description of RANSAC see Appendix~\ref{app:ransaclinreg}.}
\If{RANSAC fails}
    \State \textbf{fail}.
\EndIf
\State Fit $(a,b)$ by total least squares on $\mathcal I_\star$. \Comment{$\mathcal I_\star$ is the inlier set of RANSAC}.
\State Normalise $(a,b)$ so that $\|a\|_2=1$ and fix a deterministic sign convention.
\State Recompute the refined consensus set
\Statex\hspace{1.5em}$
\mathcal I
\gets
\{x\in\mathcal C:|a^\top x+b|\leq\varepsilon\}.
$
\If{$|\mathcal I|<m_{\rm req}$}
    \State \textbf{fail}.
\EndIf
\State \Return $(a,b)$.
\end{algorithmic}
\end{algorithm}

A simple extension for recovering multiple switching hyperplanes is to apply the procedure sequentially, removing inliers associated with the strongest hyperplane recovered before fitting the next one. This produces an arrangement of recovered affine interfaces, rather than an explicit reconstruction of the polygonal cells, and is intended for well-separated interfaces. We illustrate this sequential extension and a support-based estimate of the number of hyperplanes in Appendix~\ref{app:three-line-appendix}.

\subsection{Learning the dynamics}

Once an approximate geometry has been obtained, we need to approximate the dynamics. In this work, geometry recovery and dynamics learning are performed sequentially: after estimating the polyhedral partition $P$, we fix it during the training of the vector field. If the switching geometry is known from modelling considerations, the geometry-learning stage can be skipped, and the same dynamics-learning procedure can be applied directly on the prescribed partition. 

We denote the recovered partition by $\widetilde{P}=\{\tp_1,\ldots,\tp_{\widetilde K}\}$, allowing its number of regions $\widetilde K$ to differ from that of the target partition $K$. In our experiments, we assume we know the number of hyperplanes. However, in Appendix \ref{app:three-line-appendix}, we show that such a number can be detected during the geometry-learning phase. The goal now is to define an approximate vector field of the form
\[
\dot{x} = \tilde{f}(x,t):=\tilde{f}_k(x,t),\quad x\in\tp_k,\,\,k=1,\ldots,\widetilde K.
\]
Depending on the available structural knowledge, we choose a class $\mathcal{F}_{\widetilde{P}}^d$ for $\tilde f$. We call this class \emph{adapted to the partition $\widetilde P$} if every element has discontinuities only on the hyperplanes defining the partition $\widetilde P$. This prevents the learned vector field from creating spurious switching hyperplanes away from the recovered geometry. For example, if the vector field is known to be autonomous and piecewise affine, we could set
\[
\mathcal{F}_{\widetilde{P}}^d = \left\{\tilde{f}:\R^d \to \R^d\mid\, \tilde f(x)=\sum_{k=1}^{\widetilde K} 1_{\tp_k}(x)(B_k x+c_k),
\qquad
B_k\in\mathbb R^{d\times d},\quad c_k\in\mathbb R^d\right\}.
\]

When lacking further structural knowledge on the target dynamics, we model the local maps $\tilde f_k$ using sufficiently flexible neural networks. In \Cref{se:approxTheory}, we show how to realise such partition-adapted neural vector fields without explicitly evaluating indicator functions.

Once an approximation space is defined, we use the available data to recover a vector field that reproduces the behaviour observed in the training data. To do so, we minimise $\mathcal L_{\mathrm{traj}}(g)$ over $g\in\mathcal{F}_{\widetilde{P}}^d$, where
\begin{equation}\label{eq:lossTraj}
\mathcal L_{\mathrm{traj}}(g)
=
\frac{1}{NJd}\sum_{n=1}^N \sum_{j=1}^{J} \left\|y_j^n - \varphi_g(t_j^n;t_0^n,\xi_n) \right\|^2,
\end{equation}
$\varphi_g(t;t_0,\xi)$ denotes a numerical approximation at time $t$ of the corresponding solution of
\[
\dot z(\tau)=g(z(\tau),\tau),
\qquad
z(t_0)=\xi.
\]
To improve the training process, we increase the length of the trajectory segments $J$ gradually during training. We call this curriculum training and provide more details in \Cref{sec:experiments}. Here we set $\xi_n=y_0^n$, the noisy observed initial condition. The numerical method $\varphi$ that we consider is described in \Cref{sec:experiments}, and the precise class of neural vector fields that we use is defined in \Cref{sec:theory}. 

In some numerical experiments, we augment the trajectory-matching loss above with velocity-based auxiliary terms. Let $\mathcal I$ denote the set of trajectory indices $(n,j)$ considered by the auxiliary losses. Let $v_{h,j}^n$ denote a velocity estimate computed from the observed trajectory using one of the linear estimators in \Cref{eq:central_diff,eq:savgol,eq:lin_est}. We define
\begin{equation}\label{eq:lossVF}
\mathcal L_{\rm vf}(g)
=
\frac{1}{d|\mathcal I_{\rm vf}|}
\sum_{(n,j)\in\mathcal I_{\rm vf}}
\left\|
    g(y_j^n,t_j^n)-v_{h,j}^n
\right\|^2,
\end{equation}
where $\mathcal I_{\rm vf}\subseteq\mathcal{I}$ denotes the set of samples used for vector-field supervision and $|\mathcal I_{\rm vf}|$ the number of such samples.  This set may exclude samples close to the switching set, where velocity estimates are least reliable.

The numerical experiments using these directional regularisers have a single switching hyperplane with unit normal $a$. For arrangements of several hyperplanes, analogous penalties could be computed for each hyperplane and averaged. Let
\[
\Pi_{\mathrm{tan}} = I_d-aa^\top
\]
be the orthogonal projection onto the tangent space of the hyperplane, and define
\[
\mathcal I_{\mathrm{nor}}
=
\{(n,j)\in\mathcal I: |v_{h,j}^n\cdot a|>\tau_{\mathrm{nor}}\},
\qquad
\mathcal I_{\mathrm{tan}}
=
\{(n,j)\in\mathcal I: \|\Pi_{\mathrm{tan}}\, v_{h,j}^n\|_2>\tau_{\mathrm{tan}}\}.
\]
The normal regulariser penalises a predicted normal component that points in the opposite direction to the estimated velocity,
\begin{equation}\label{eq:lossNor}
\mathcal L_{\mathrm{nor}}(g)
=
\frac{1}{|\mathcal I_{\mathrm{nor}}|}
\sum_{(n,j)\in\mathcal I_{\mathrm{nor}}}
\max\left\{
0,\,
-\bigl(g(y_j^n,t_j^n)\cdot a\bigr)
\operatorname{sign}(v_{h,j}^n\cdot a)
\right\},
\end{equation}
while the tangential regulariser matches the tangential components,
\begin{equation}\label{eq:lossTan}
\mathcal L_{\mathrm{tan}}(g)
=
\frac{1}{d|\mathcal I_{\mathrm{tan}}|}
\sum_{(n,j)\in\mathcal I_{\mathrm{tan}}}
\left\|
    \Pi_{\mathrm{tan}} \,g(y_j^n,t_j^n)
    -
    \Pi_{\mathrm{tan}} \,v_{h,j}^n
\right\|^2.
\end{equation}
The general objective used in the numerical experiments is therefore
\begin{equation}\label{eq:lossFull}
\mathcal J(g)
=
\lambda_{\mathrm{traj}}\mathcal L_{\mathrm{traj}}(g)
+
\lambda_{\mathrm{vf}}\mathcal L_{\mathrm{vf}}(g)
+
\lambda_{\mathrm{nor}}\mathcal L_{\mathrm{nor}}(g)
+
\lambda_{\mathrm{tan}}\mathcal L_{\mathrm{tan}}(g),
\end{equation}
with the weights, thresholds, and active terms specified in \Cref{sec:experiments}.

\section{Theoretical analysis}\label{sec:theory}
We now analyse certain theoretical aspects of the methodology. First, we study how observational noise interplays with the properties of jumps across the switching hyperplane in the context of geometry learning, and then we study the approximation properties of a class of piecewise-smooth neural networks later used in the experiments. 

\subsection{Analysis of geometry learning}\label{sec:theory_geometryLearning}
In this section, we address some theoretical aspects of the geometry learning described in \Cref{sec:geom_learning}. Here, we consider the simplified case where $\R^d$ is separated by one hyperplane, i.e., we suppose that 
\[
  \Sigma:=\left\{x\in\mathbb{R}^d:\ s(x)=0\right\},
  \qquad s(x):=a^\top x+b,
\]
where $a\in\mathbb{R}^d\setminus\left\{0\right\}$ and $b\in\mathbb{R}$, and the goal is to recover $a$ and $b$. 

Velocity-based discontinuity learning with linear velocity estimation methods may be sensitive to noise. Since the observations are affected by random noise, they are random variables, and consequently so are the velocity estimates. As linear velocity estimators rely on numerical differentiation, they amplify high-frequency components of the data, such as noise.
The detection of $\Sigma$ also relies on locally amplifying the velocity jump. 
If the noise magnitude is too high, a linear velocity estimator cannot detect jumps across the boundary.
The amplification has a floor that no choice of weights can beat.
For an estimator $v_{h,j}^n = \sum_{k=-r}^r w_k y_{j+k}^n$ applied to data of the form \eqref{eq:data}, the random part $\sum_{k=-r}^r w_k \eta_{j+k}^n$ is a linear combination of independent Gaussians, so
\[
    v_{h,j}^n - \sum_{k=-r}^r w_k x_*^n(t_{j+k}^n) \sim N\left(0,\sigma^2 \|w\|_2^2 I_d\right), \qquad \|w\|_2^2 = \sum_{k=-r}^r w_k^2 .
\]
The consistency conditions $\sum_k w_k = 0$ and $\sum_k w_k kh = 1$ require, by Cauchy--Schwarz,
\[
    1 = \left|\sum_{k=-r}^r w_k kh\right| \leq h \sqrt{\sum_{k=-r}^r w_k^2}\,\sqrt{\sum_{k=-r}^r k^2} = h\|w\|_2\sqrt{S_r}, \qquad S_r = \frac{r(r+1)(2r+1)}{3} ,
\]
so $\|w\|_2 \geq 1/(h\sqrt{S_r})$. This floor decreases as the window radius $r$ grows, but increases as the step size $h$ shrinks. It explains why a Savitzky--Golay estimator is more robust to noise than central differences, and also why it will still struggle when $\sigma$ is too large.

The noise floor on a single velocity estimate has consequences for which jumps across $\Sigma$ are detectable. 
The next result combines this floor with the smoothing error incurred by the choice of window size $r$, addressing which jump sizes are detectable. 
To this end, consider the estimated velocity jump
\begin{align*}
    \Delta v_{h,j}^n = v_{h,j+1}^n - v_{h,j}^n = \sum_{k =-r}^r w_k y^n_{j+1+k} -\sum_{k =-r}^r w_k y^n_{j+k} = \sum_{k =-r}^{r+1} c_k y^n_{j+k}, 
\end{align*}
where 
\begin{align*}
    c_k = \begin{cases}
        -w_{-r} \quad \text{if } k = -r,\\
        w_{r} \quad \text{if } k = r+1, \\
        w_{k-1}-w_{k} \quad \text{else}.
    \end{cases}
\end{align*}
If the velocity estimator satisfies the same consistency conditions $\sum_k w_k = 0$, $\sum_k w_k kh = 1$, then $\sum_{k=-r}^{r+1} c_k = 0$ and $\sum_{k=-r}^{r+1} c_k kh = 0$, i.e., no jump is detected on constant and constant-velocity trajectories. Further, it holds that $\|c\|_2 \leq 2\|w\|_2$.
 We are now ready to state the proposition. 
\begin{proposition}
    Consider the velocity estimator $v_{h,j}^n = \sum_{k =-r}^r w_k y^n_{j+k}$ satisfying the consistency conditions $\sum_k w_k = 0$, $\sum_k w_k kh = 1$. Let $f^+$ and $f^-$ be smooth and extendable to the whole domain. Assume that a reference trajectory $x^n_*$ crosses $\Sigma$ transversally at one single time $\tau$, and that without loss of generality, $s(x)<0$ before the crossing. Set $p = x^n_*(\tau)$$\in \Sigma$. Denote by $x^+$ the solution to $\dot x = f^+(x,t)$, with initial time $\tau$ and starting value $x^+(\tau) = p$, and denote by $x^-$ the trajectory of $f^-$ agreeing with $x^n_*$ before the switch. By continuity, $x^-(\tau) = x^+(\tau) = p$. The velocity jump is given by
    \begin{align*}
        f^+(p,\tau)-f^-(p,\tau) = \kappa e, \quad \kappa>0, \; e\in \mathbb{R}^d,
    \end{align*}
    where $\kappa = \|f^+(p,\tau)-f^-(p,\tau)\|_2$ and $\|e\|_2=1$. Assume that the second derivatives of both $x^+$ and $x^-$ are uniformly bounded by $M$. Let
    \begin{align*}
        &R =  \sum_{k= -r}^{r+1} c_k \max\left(t_{j+k}-\tau,0 \right), ~\alpha = M \sum_{k:t_{j+k}>\tau} |c_k| (t_{j+k}-\tau)^2, ~ \beta =  \frac{M}{2} \sum_{k=-r}^{r+1} |c_k| (t_{j+k}-t_j)^2 .
    \end{align*}
    Consider now two scenarios. In the first, the trajectory jumps in the stencil used to compute $\Delta v_{h,j}^n$, and in the second, it follows $x^-$ throughout. For simplicity, denote the first jump estimate by $\Delta_1$, and the second by $\Delta_2$. Then,
    \begin{enumerate}
        \item Both $\Delta_1$ and $\Delta_2$ are Gaussian, that is, $\Delta_i \sim \mathsf{N}(\mu_i,C)$, with
        \[C =  \sigma^2 \|c\|_{2}^2 I_d,~ \mu_1-\mu_2 = \kappa R e + \varrho,\] where $\varrho \in \mathbb{R}^d$ satisfies $\|\varrho\|_2 \leq \alpha$.
        \item The expected magnitudes satisfy
        \begin{align*}
            \mathbb{E}[\|\Delta_1\|_2] \geq \kappa|R|-\beta-\alpha, \quad \mathbb{E}[\|\Delta_2\|_2] \leq \beta + \sigma \|c\|_{2} \sqrt{d},
        \end{align*}
        so the expected jump magnitude exceeds the expected smooth magnitude whenever $\kappa |R| > 2\beta +\alpha+\sigma \|c\|_{2} \sqrt{d}$.
    \end{enumerate}
\end{proposition}
The constant $R$ measures how much of the post-crossing displacement the jump stencil captures.
The constants $\alpha$ and $\beta$ are second-derivative Taylor remainders, both controlled by $M$: $\beta$ bounds the spurious jump produced by a smooth trajectory of nonzero acceleration, while $\alpha$ bounds the error in approximating the post-crossing displacement by its leading linear-in-time term.
Detectability therefore requires the scaled jump $\kappa|R|$ to dominate $\alpha$, $\beta$, and the noise term $\sigma\|c\|_2\sqrt{d}$.
\begin{proof}
    As above, the randomness is introduced by the estimator being applied to the observational noise. Due to the linearity, the resulting quantity is also normally distributed, that is,
    \[
    \sum_{k=-r}^{r+1} c_k \eta_{j+k}^n \sim \mathsf{N}(0,\sigma^2 \|c\|_{2}^2 I_d).
    \]
    Consider now $\Delta_2$, i.e., when the stencil contains no jumps. In that case, we Taylor expand and write
    \[
    x^-(t_{j+k}) = x^-(t_j) + kh\dot{x}^-(t_j) + r_k,
    \]
    where $r_k = \int_{t_j}^{t_{j+k}} (t_{j+k}-s)\ddot{x}^-(s)\,\mathrm{d}s$, so that $\|r_k\|_2$ $\leq \tfrac{M}{2}(t_{j+k}-t_j)^2$. The annihilation property of the velocity jump estimator means that
    \[
    \mu_2 = \sum_{k=-r}^{r+1} c_k r_k,
    \]
    which satisfies
    \[
    \|\mu_2\|_2 \leq \frac{M}{2} \sum_{k= -r}^{r+1} |c_k| (t_{j+k}-t_j)^2 = \beta.
    \]
    The case when there is a jump is similar. The deterministic trajectory is $x^-(t) + (x^+(t)-x^-(t))\mathbf{1}_{t > \tau}$, where $\mathbf{1}_{t > \tau}$ is $1$ if $t > \tau$ and $0$ elsewhere. Since $x^+(t)-x^-(t)$ is smooth, we Taylor expand around $\tau$ and obtain
    \[
    x^+(t_{j+k})-x^-(t_{j+k}) = \kappa e (t_{j+k}-\tau) + \varrho_k,
    \]
    where $\varrho_k = \int_\tau^{t_{j+k}} (\ddot{x}^+(s)-\ddot{x}^-(s))(t_{j+k}-s)\,\mathrm{d}s$. Since $\|\ddot{x}^+ - \ddot{x}^-\|_2$ $\leq 2M$, we see that $\|\varrho_k\|_2$ $\leq M(t_{j+k}-\tau)^2$. Summing this up, we obtain
    \[
    \mu_1 = \mu_2 + \sum_{k:t_{j+k}>\tau} c_k\left[\kappa e (t_{j+k}-\tau) + \varrho_k\right] = \kappa e R + \sum_{k:t_{j+k}>\tau} c_k\varrho_k,
    \]
    where we have that
    \[
    \left\|\sum_{k:t_{j+k}>\tau} c_k\varrho_k\right\|_2 \leq M \sum_{k:t_{j+k}>\tau} |c_k| (t_{j+k}-\tau)^2 = \alpha.
    \]
    For the second statement, first note that by the triangle inequality and Jensen's inequality,
    \begin{align*}
        \mathbb{E}[
        \|\Delta_2\|_2] \leq \|\mu_2\|_2 + \sigma \|c\|_{2}\sqrt{d} \leq \beta + \sigma \|c\|_{2}\sqrt{d}.
    \end{align*}
    For the second bound, Jensen's inequality using the norm function, followed by the reverse triangle inequality, gives
    \begin{align*}
        \mathbb{E}[
        \|\Delta_1\|_2] \geq \|\mathbb{E}[\Delta_1]\|_2 \geq \kappa |R| - \beta - \alpha.
    \end{align*}
\end{proof}
In practice, this means that the gap between the jump magnitude and the noise standard deviation must be sufficiently large to distinguish smooth points from true jumps, and that sensitivity decreases as the dimension grows.
We numerically test the balance between the discontinuity gap and the noise magnitude in \Cref{sec:experiments}.

\begin{remark}
The proposition above is a local signal-detection result.  It assumes that the crossing time $\tau$ lies in the range of the velocity-jump stencil and gives a condition under which the resulting derivative jump is distinguishable from smooth curvature and observation noise.  The subsequent hyperplane-fitting step uses the corresponding neighbouring trajectory points as candidates for $\Sigma$.  This introduces a separate conditioning issue. Points with small signed distance to $\Sigma$ need not be tightly localised around the crossing time if the normal crossing velocity is small. 

Indeed, let $p=x(\tau)$ and suppose, for simplicity, that $\Sigma=\{x:s(x)=0\}$ is a hyperplane with $s(x)=a^\top x+b$ and $\|a\|_2=1$. For a nearby time $t$ on either side of the crossing,
\[
s(x(t)) = (t-\tau)a^\top f^{\rm pre/post}(p,\tau) + \mathcal O(|t-\tau|^2).
\]
Thus, the signed-distance signal scales like $|t-\tau|v_{\rm nor}$, where
\[
v_{\rm nor}=\min\left\{|a^\top f^{\rm pre}(p,\tau)|,\,|a^\top f^{\rm post}(p,\tau)|\right\}.
\]
With observations $y(t)=x(t)+\eta(t)$ and $\eta(t)\sim\mathsf N(0,\sigma^2I_d)$, one has
\[
s(y(t))-s(x(t))=a^\top\eta(t)\sim\mathsf N(0,\sigma^2).
\]
Consequently, the time interval over which the clean signed distance is comparable to the observation noise has a width of order $\sigma/v_{\rm nor}$. Indeed, this is the regime where $|t-\tau|v_{\rm nor}\approx\sigma$, so the observed signed distance $s(y(t))$ is strongly affected by noise and is no longer a reliable proxy for the clean distance to $\Sigma$. When $v_{\rm nor}$ is small, a larger portion of the trajectory near the crossing is ambiguous from the point of view of geometric localisation.  Thus, the jump-to-noise condition above should be read together with a non-grazing condition: the normal crossing velocity must be large enough, relative to the sampling step and observation noise, for the candidate cloud to localise $\Sigma$ reliably.
\end{remark}

\subsection{Universal neural network for piecewise-smooth dynamical systems} \label{se:approxTheory}

Let us fix a polyhedral partition $P=\{\p_1,\dots,\p_K\}$ of $\R^d$ induced by a finite affine hyperplane arrangement. We now define a neural class $\mathcal F_P^d$ adapted to this partition. The class has two key properties: its discontinuities are confined to the prescribed arrangement, and it can approximate piecewise-smooth vector fields arbitrarily well almost everywhere on compact sets.

Let $\operatorname{ReLU}(x)=\max\{0,x\}$ denote the ReLU function, and $\mathcal{H}(x)$ the Heaviside step function. We use the convention $\mathcal H(r)=1$ for $r\ge 0$ and $\mathcal H(r)=0$ for $r<0$, applied componentwise to vectors. For a vector $J\in\R^k$, we define the nonlinear function $\operatorname{ReLU}_J:\R^k\to\R^k$ by $\operatorname{ReLU}_J(x)=\operatorname{ReLU}(x)+J\odot \mathcal{H}(x)$, as in \cite{della2023discontinuous}. For fixed $d_1,d_2\in\mathbb{N}$, we introduce the networks
\begin{align*}
\mathcal{F}_{\operatorname{ReLU}_J}^{d_1,d_2}:=\Bigl\{\mathcal N:\R^{d_1}\to\R^{d_2}\,\Big|\;&\mathcal N=W_L\circ \operatorname{ReLU}_{J_L}\circ \cdots \circ W_1\circ \operatorname{ReLU}_{J_1}\circ W_0,\\
&L\in\mathbb N,W_\ell:\R^{h_\ell}\to\R^{h_{\ell+1}}\text{ affine},\,\,\ell=0,\dots,L\\
&J_{\ell}\in\R^{h_{\ell}},\,\,\ell=1,\dots,L,\,\,h_0=d_1,\,\, h_{L+1}=d_2,\,\, h_1,\dots,h_L\in\mathbb N\Bigr\}.\nonumber
\end{align*}
We denote by $\mathcal F_{\operatorname{ReLU}}^{d_1,d_2}\subset \mathcal F_{\operatorname{ReLU}_J}^{d_1,d_2}$ the subclass for which all jump parameters vanish, i.e., the set of deep $\mathrm{ReLU}$ feedforward networks. In this section, we repeatedly use the identities
\begin{equation}\label{eq:heaviJump-note}
\mathcal H(x)=
\begin{bmatrix}I_d&-I_d\end{bmatrix}
\operatorname{ReLU}_J \left(
\begin{bmatrix}I_d\\ I_d\end{bmatrix}x
\right),
\qquad x\in\R^d,\,
J=
\begin{bmatrix}1_d\\ 0_d\end{bmatrix}\in\R^{2d},
\end{equation}
and
\begin{equation}\label{eq:xJump-note}
x=
\begin{bmatrix}I_d&-I_d\end{bmatrix}
\operatorname{ReLU}_J\left(
\begin{bmatrix}I_d\\ -I_d\end{bmatrix}x
\right),
\qquad x\in\R^d,\,
J=
\begin{bmatrix}0_d\\ 0_d\end{bmatrix}\in\R^{2d}.
\end{equation}
In particular, $\mathcal H\in \mathcal F_{\operatorname{ReLU}_J}^{d,d}$ and the identity map belongs to $\mathcal F_{\operatorname{ReLU}}^{d,d}$.
\begin{lemma}[Closure under affine combinations and parallel concatenation]\label{lem:closure-basic}
For every choice of dimensions, both $\mathcal F_{\operatorname{ReLU}_J}^{d_1,d_2}$ and $\mathcal F_{\operatorname{ReLU}}^{d_1,d_2}$ are closed under multiplication with a scalar, linear combination, and parallel concatenation.
\end{lemma}
The proof follows similar ideas as in \cite[Section 5.1]{petersen2024mathematical}, and it can be found in Appendix \ref{app:additionalProofsUniversal}.

Adopting the notation in \eqref{eq:partitioning}, we set
\[
d_P:=\sum_{k=1}^K m_k,
\qquad
A_P:=\begin{bmatrix}A_1\\ \vdots \\ A_K\end{bmatrix},
\qquad
b_P:=\begin{bmatrix}b_1\\ \vdots \\ b_K\end{bmatrix},
\qquad
T_P(x):=b_P-A_Px\in\R^{d_P},
\]
and define $z_P(x)=\mathcal{H}(T_P(x))\in\{0,1\}^{d_P}$, 
\[
\Gamma_P
:=
\bigcup_{\ell=1}^{d_P}
\left\{x\in\R^d:\ (T_P(x))_\ell=0\right\}.
\]
Thus $z_P$ is locally constant on $\R^d\setminus\Gamma_P$.

\begin{definition}[Geometry-preserving class associated with $P$]\label{def:FP}
For the fixed state dimension $d$, we define $\mathcal F_P^d$ as the set of maps
\[
(x,t)\mapsto \Psi\bigl(x,t,\Phi(z_P(x))\bigr)
\]
such that
\[
\Phi\in \mathcal F_{\operatorname{ReLU}_J}^{d_P,d_1},
\qquad
\Psi\in 
\mathcal F_{\operatorname{ReLU}}^{1+d+d_1,d},
\qquad
d_1\in\mathbb N.
\]
\end{definition}
\begin{remark}
Time-independent maps of the form $x\mapsto \Psi(x,\Phi(z_P(x)))$ with $\Psi\in \mathcal{F}_{\operatorname{ReLU}}^{d+d_1,d}$, $\Phi\in \mathcal{F}_{\operatorname{ReLU}_J}^{d_P,d_1}$, $d_1\in\mathbb{N}$, can also be realised as maps in $\mathcal{F}_P^d$. If $\mathcal{N}\in\mathcal{F}_P^d$ is time-independent, we write $\mathcal{N}(x)$ and avoid including the variable $t$ for conciseness.
\end{remark}
\begin{proposition}[No new jumps]\label{prop:no-new-jumps}
Every element of $\mathcal F_P^d$ is continuous on $(\R^d\setminus \Gamma_P)\times\R$.
\end{proposition}
\begin{proof}
Let $\mathcal N\in\mathcal{F}_P^d$. If $x_0\notin \Gamma_P$, then $z_P$ is locally constant near $x_0$. Therefore, near $x_0$, the map $x\mapsto \Phi(z_P(x))$ is constant, while $(x,t)\mapsto \Psi(x,t,\Phi(z_P(x)))$ is the composition of continuous maps. Thus $\mathcal N$ is continuous at $(x_0,t)$ for every $t\in\R$.
\end{proof}

\begin{proposition}[Linearity of the geometry-preserving class]\label{prop:FP-linear}
$\mathcal F_P^d$ is a vector space over $\R$.
\end{proposition}
The proof follows from the vector space property of $\mathcal{F}_{\operatorname{ReLU}_J}^{d_1,d_2}$. We include it in Appendix \ref{app:additionalProofsUniversal}.

\begin{theorem}[Exact ReLU gating on compact sets]\label{thm:reluGating}
Let $\Omega\subset\mathbb R^d$ be compact, $T>0$, and $R\in \mathcal{F}_{\operatorname{ReLU}}^{d+1,m}$. Then there exists $\widetilde R\in \mathcal{F}_{\operatorname{ReLU}}^{d+2,m}$ such that
\[
\widetilde R(x,t,s)=
\begin{cases}
R(x,t),& s=1,\\
0,& s=0,
\end{cases}
\qquad \forall\, (x,t,s)\in \Omega\times[0,T]\times\{0,1\}.
\]
\end{theorem}
Since affine maps of the form $(x,t)\mapsto Ax+b\in\R^m$ belong to $\mathcal{F}_{\operatorname{ReLU}}^{d+1,m}$, the result above applies to affine time-independent maps as well.

\begin{proof}
Since $R$ is continuous, and $[0,T]$ and $\Omega$ are compact, \[
\lambda:=\sup_{\substack{x\in\Omega\\ t\in [0,T]}}\|R(x,t)\|_\infty<\infty.
\]
We propagate the binary variable $s$ unchanged through the network while computing $R(x,t)$ in parallel. Since $s\in\{0,1\}\subset [0,\infty)$, one has $\operatorname{ReLU}(s)=s$, so this can be done using only $\mathrm{ReLU}$ layers. At the final stage, define
\[
\widetilde R(x,t,s)
:=
\operatorname{ReLU}\bigl(R(x,t)-\lambda (1-s) 1_m\bigr)
-
\operatorname{ReLU}\bigl(-R(x,t)-\lambda (1-s) 1_m\bigr).
\]
If $s=1$, then
\[
\widetilde R(x,t,1)=\operatorname{ReLU}\bigl(R(x,t)\bigr)-\operatorname{ReLU}\bigl(-R(x,t)\bigr)=R(x,t).
\]
If $s=0$, then for every component one has $|R_j(x,t)|\le \lambda$, hence $\widetilde R(x,t,0)=0$ as desired.
\end{proof}

We now show that the networks in $\mathcal{F}_P^d$ can approximate  every piecewise-smooth vector field within an arbitrary accuracy.

\begin{definition}[Piecewise-smooth vector field on a fixed arrangement]
Let $\Omega\subset\R^d$ be compact and let $P=\{\p_1,\dots,\p_K\}$ be as above. A function $F:\Omega\times\R\to\R^d$ is said to be piecewise-smooth with respect to $P$ if for every $k\in\{1,\dots,K\}$ there exists a map
\[
F_k:\overline{\Omega\cap \p_k}\times\R\to\R^d
\]
which is continuous on $\overline{\Omega\cap \p_k}\times\R$ and smooth on $(\Omega\cap \p_k)\times\R$, and such that
\[
F(x,t)=F_k(x,t),
\qquad \forall\,x\in\Omega\cap \p_k,\,t\in\R.
\]
\end{definition}
\begin{theorem}[Approximation of piecewise-smooth systems without creating new jumps]\label{thm:approx-piecewise-smooth}
Let $\Omega\subset\R^d$ be compact, $T>0$, and $P=\{\p_1,\dots,\p_K\}$ a finite family of pairwise disjoint open convex polyhedra with
\[
\Omega\setminus \Gamma_P = \bigcup_{k=1}^K \left(\Omega\cap \p_k\right).
\]
Let $F:\Omega\times\R\to\R^d$ be piecewise smooth with respect to $P$. Then for every $\varepsilon>0$, there exists $\mathcal N\in\mathcal F_P^{d}$ such that
\[
\sup_{\substack{x\in \Omega\setminus\Gamma_P \\ t\in [0,T]}}\left\|F(x,t)-\mathcal N(x,t)\right\|_\infty<\varepsilon.
\]
\end{theorem}
\begin{proof}
For every $k=1,\ldots,K$, we have that $F_k:\overline{\Omega\cap\p_k}\times[0,T]\to\R^d$ is continuous. Since $\overline{\Omega\cap\p_k}\times[0,T]$ is compact, the universal approximation theory of $\mathrm{ReLU}$ networks, see \cite{yarotsky2017error}, ensures that there exists $\widetilde\Psi_k \in \mathcal{F}_{\operatorname{ReLU}}^{d+1,d}$ such that
\[
\sup_{\substack{x\in \overline{\Omega\cap\p_k} \\ t\in [0,T]}}\left\|F_k(x,t) - \widetilde \Psi_k(x,t)\right\|_\infty<\varepsilon.
\]
Consider the projection $\Pi_k:\R^{d_P}\to\R^{m_k}$ onto the block corresponding to $\p_k$, so that
\[
\Pi_k T_P(x)=b_k-A_k x.
\]
We define the function $\chi_k(z)=\mathcal H\!\left(1_{m_k}^\top \Pi_k z-m_k+\frac12\right)$ in $\mathcal{F}_{\operatorname{ReLU}_J}^{d_P,1}$, compare \cite[Equation 1]{kong2025expressivitydeepheavisidenetworks}. If $x\notin\Gamma_P$, then $x\in\p_k$ if and only if all entries of $\Pi_k z_P(x)=\mathcal H(b_k-A_kx)$ are equal to $1$. Therefore, $1_{\p_k}(x) = \chi_k(z_P(x))$ for $x\in\Omega\setminus\Gamma_P$. 

By \Cref{thm:reluGating}, there exists a ReLU network $\Psi_k\in\mathcal F_{\operatorname{ReLU}}^{d+2,d}$ such that
\[
\Psi_k(x,t,s)=
\begin{cases}
\widetilde\Psi_k(x,t), & s=1,\\
0, & s=0,
\end{cases}
\qquad \forall\,(x,t,s)\in \Omega\times[0,T]\times \{0,1\}.
\]
Therefore the map
\[
\mathcal N_k(x,t):=\Psi_k\bigl(x,t,\chi_k(z_P(x))\bigr)
\]
belongs to $\mathcal F_P^d$ and satisfies
\[
\mathcal N_k(x,t)=1_{\p_k}(x)\widetilde\Psi_k(x,t),
\qquad \forall\,(x,t)\in (\Omega\setminus\Gamma_P)\times[0,T].
\]
Since $\mathcal F_P^d$ is a vector space, the sum
\[
\mathcal N(x,t):=\sum_{k=1}^K \mathcal N_k(x,t)
\]
belongs to $\mathcal F_P^d$. Moreover, for $(x,t)\in(\Omega\setminus\Gamma_P)\times[0,T]$,
\[
\mathcal N(x,t)
=
\sum_{k=1}^K 1_{\p_k}(x)\widetilde\Psi_k(x,t).
\]
For each $x\in\Omega\setminus\Gamma_P$, there is a unique $\underline{k}\in\{1,\ldots,K\}$ such that $x\in\p_{\underline{k}}$. Therefore, 
\[
\left\|\mathcal{N}(x,t) - F(x,t)\right\|_\infty = \left\|\widetilde \Psi_{\underline{k}}(x,t)-F(x,t)
\right\|_\infty < \varepsilon
\]
as desired.
\end{proof}

\begin{remark}[Exact representation of piecewise-affine vector fields on compacts]
Given a partition $P=\{\p_1,\dots,\p_K\}$ and a compact set $\Omega\subset\R^d$ as in \Cref{thm:approx-piecewise-smooth}, every piecewise-affine function
\[
x\mapsto \sum_{k=1}^K 1_{\p_k}(x)(W_kx+w_k),\,\,W_1,\dots,W_K\in\R^{d\times d},\,\,w_1,\dots,w_K\in\R^d
\]
is the restriction over $\Omega$ of a function $\mathcal{N}\in\mathcal{F}_P^d$. Such a result follows by picking $\widetilde\Psi_k(x,t) = W_kx+w_k$ in the proof above.
\end{remark}

\begin{remark}[Extended approximation spaces]
Replacing the raw time input by a continuous feature map \(\varphi(t)\), or augmenting \(\Psi\) with other continuous features of \(x\) and \(t\), does not enlarge the possible discontinuity set because \(z_P(x)\) remains the only discontinuous input. In the
periodically forced experiments, we therefore use maps of the form
\[
(x,t)\mapsto
\Psi\bigl(x,\varphi(t),\Phi(z_P(x))\bigr),
\]
where
\[
\Psi\in \mathcal F_{\operatorname{ReLU}}^{d+2k+d_1,d},
\qquad
\Phi\in \mathcal F_{\operatorname{ReLU}_J}^{d_P,d_1},
\]
and
\[
\varphi(t)
=
\begin{bmatrix}
\cos(\omega_1 t + \theta_1)&
\sin(\omega_1 t + \theta_1)&
\dots &
\cos(\omega_k t + \theta_k)&
\sin(\omega_k t + \theta_k)
\end{bmatrix}^\top
\in\R^{2k}.
\]
The frequencies $\omega_1,\dots,\omega_k>0$ may be fixed from prior knowledge
or chosen as random Fourier features.
\end{remark}

\section{Numerical experiments}\label{sec:experiments}

We now test the methodology introduced and analysed above. The experiments
address three questions: (1) whether the proposed architecture confines
discontinuities to a prescribed partition, (2) whether derivative-jump RANSAC can
recover the switching geometry from the training trajectories, and (3) whether the recovered
geometry supports accurate trajectory learning. We study the geometry and
dynamics learning phases both in isolation, in \Cref{sec:regression}
and \Cref{sec:learn-geometry}, and when combined, in
\Cref{sec:full_pipeline}. We also test the effect of incorporating additional
prior knowledge, such as the second-order structure of the target dynamics and
the forcing frequency in periodically forced systems. Unless stated otherwise, all experiments based on training a network use the random seed $7$. 

\vspace{0.1in}

\noindent In the experiments in which we recover the geometry from trajectories, we focus on systems with one
dominant affine switching surface. The dry-friction and the climate-based PP04 benchmarks both have a physically meaningful single switch, and this
setting isolates the interaction between derivative-jump recovery and dynamics
learning without adding the separate task of assigning several recovered facets
to regions. The approximation architecture is formulated for general finite
polyhedral partitions, and the known-geometry regression experiment already
tests a multi-region partition. The more general case is discussed in Appendix~\ref{app:three-line-appendix}, which gives a
synthetic demonstration of sequential RANSAC--TLS recovery for multiple affine
interfaces.

The GitHub repository associated with the paper is at \url{https://github.com/davidemurari/learn-piecewise-smooth}.

\vspace{0.1in}

\paragraph{Benchmark dynamical systems}
The three benchmark systems are chosen to test different parts of the methodology.

\vspace{0.1in}

\noindent {\bf Benchmark A:} We first isolate the dynamics-learning architecture on a two-dimensional
piecewise-linear oscillator with prescribed switching geometry. The vector field is
\begin{equation}\label{eq:pw-oscillator}
    F(x)=
    \begin{bmatrix}
        p\\
        -q-\alpha(q)p
    \end{bmatrix},
    \qquad
    \alpha(q)=
    \begin{cases}
        0.9, & q<-1,\\
        -0.9, & -1\le q<1,\\
        0.9, & q\ge 1.
    \end{cases}
\end{equation}
This example has three smooth regions separated by the two switching lines
$q=\pm1$, and is used to test whether the proposed architecture can approximate
a discontinuous vector field without introducing spurious jumps away from the
prescribed partition.

\vspace{0.1in}

\noindent We then test the complete learning pipeline on two families of trajectory-learning
problems: dry-friction oscillators and the PP04 system.

\vspace{0.1in}

\noindent {\bf Benchmark  B: The dry-friction
oscillator.} This system, also considered in \cite{lathourakis2024physics} and analysed in \cite{brogliato2016nonsmooth,Laurea2008book}, models the stick-slip dynamics of a forced mass moving on a surface with friction, and is useful because
it combines a simple switching geometry with non-trivial Filippov behaviour. It
allows us to separate three effects: the benefit of knowing the switching
geometry, the benefit of imposing the mechanical second-order structure, and the
degradation caused by noisy trajectory observations. The system has state
$x=(q,p)$ and stick-slip dynamics given by
\begin{equation}\label{eq:dry-friction}
\dot q = p,
\qquad
\dot p = -q - c p + A\cos(\omega t) - \mu\,\mathrm{sign}(p),
\end{equation}
with switching hyperplane $\Sigma=\{(q,p):p=0\}$. We consider two forced cases
with $c=0.1$, $\mu=0.5$, and $A=1$, but different forcing frequencies as in
\cite{lathourakis2024physics}:
\[
\omega=0.3 \quad \text{(case 1)}, \qquad
\omega=0.15 \quad \text{(case 2)}.
\]

\vspace{0.1in}

\noindent {\bf Benchmark C: The PP04 system.} This climate-based system models the periodic behaviour of the ice ages \cite{Paillard2004,Morupisi2020} and provides a complementary three-dimensional benchmark. In contrast to the dry-friction oscillator, its switching surface is an affine plane not aligned with a coordinate axis, and the system is periodically forced.
The state is $x=(V,A,C)$ representing the global ice volume $V$, the Antarctic ice cover $A$ and the atmospheric carbon dioxide level $C$. In this model, carbon dioxide is gradually absorbed by the oceans over a period of about 100k years, leading to global cooling and a growth of the ice sheets. At a critical concentration, the oceans can no longer absorb the carbon dioxide, and it is suddenly released into the atmosphere. This causes a sudden rise in $C$ and rapid global warming, leading to a drop in $V$ and $A$. The process then repeats, synchronised by the periodic variations in the solar insolation caused by the Milankovitch cycles.  With the parameter values used in the experiments, the linear switching function (representing the oceanic stratification linked to the carbon dioxide absorption) is given by:
\[
    s(V,A,C)=aV-bA+d,
    \qquad
    a=0.3,\quad b=0.7,\quad d=0.27.
\]
The switching surface (corresponding to the critical values of the parameters when the carbon dioxide is suddenly released from the oceans) is the plane in the 3-dimensional phase-space given by:
\[
    \Sigma=\{(V,A,C):s(V,A,C) \equiv (a,-b,0)^T(V,A,C) + d =0\}.
\]
\Cref{fig:pp04-orbits} visualises this plane in the $(V,A,C)$ phase space
and compares representative trajectories of the reference PP04 system with learned-model rollouts from the same initial conditions. Setting the (assumed periodic) solar insolation to $I_{65}(t)=\mu\sin(\omega t)$ (which varies according to the Milankovitch cycles \cite{Morupisi2020}), the vector field (where time $t$ is expressed in kilo-years) is
\begin{align}
    \dot V &=
    \frac{-xC-yI_{65}(t)+z-V}{\tau_V},\nonumber\\
    \dot A &=
    \frac{V-A}{\tau_A},\label{eq:pp04}\\
    \dot C &=
    \frac{\alpha I_{65}(t)-\beta V+\gamma H(-s(V,A,C))+\delta-C}{\tau_C}.\nonumber
\end{align}
Here the Heaviside function is defined by $H(z)=1$ for $z\geq 0$ and $H(z)=0$ otherwise. The $\gamma$ term is active when $s(V,A,C)\leq0$, including on the switching surface. In the trajectory-learning experiments all PP04 coefficients are fixed at $x=1.3$, $y=0.5$, $z=0.8$, $\alpha=0.15$, $\beta=0.5$, $\gamma=0.7$, $\delta=0.4$, $\tau_V=15$, $\tau_A=12$, $\tau_C=5$, $\mu=0.467$, and $\omega=0.1532$.  The geometry-recovery diagnostic below instead varies $\gamma$ while keeping every other coefficient and the switching plane fixed. This example tests whether the derivative-jump geometry-learning stage can accurately recover a non-axis-aligned switching plane for subsequent trajectory learning.

\subsection{Evaluation metrics}
For the different experiments, we isolate different aspects of the learning process. This is done by changing the architecture and its training regimen, as well as considering different evaluation metrics, which we introduce below.

\vspace{0.1in}

\paragraph{Discontinuity detection diagnostic}
We define a metric to detect discontinuities and test it on the system in \Cref{eq:pw-oscillator}. We describe it for this system, but it naturally extends. 

For a learned vector field $\widehat F$, output component $k$, grid point $x$, we compute
\[
     D_{h,k}(x)
    =
    \max_{v\in V}
    \left|
    \widehat F_k(x+hv)-\widehat F_k(x-hv)
    \right|,\,\,V=\left\{
    (1,0),\ (0,1),\ \frac{(1,1)}{\sqrt2},\
    \frac{(1,-1)}{\sqrt2}
    \right\}.
\]
For a continuous piecewise-linear function, $D_{h,k}(x)$ scales approximately linearly with $h$, so $D_{2h,k}(x)/D_{h,k}(x)\approx 2$. We retain a persistent jump score $P_{h,k}(x)=D_{h,k}(x)$ only where
\[
    \frac{D_{2h,k}(x)}{D_{h,k}(x)+10^{-12}} < 1.35
\]
and $D_{h,k}(x)$ exceeds a small amplitude floor. The floor is chosen componentwise as
\[
    \max\left\{10^{-6},\ 10^{-2}
    \operatorname{percentile}_{99}(D_{h,k})\right\}.
\]
We then define the known-switch band
\[
    B_\delta
    =
    \{ |q+1|\le \delta\}\cup\{|q-1|\le \delta\},
    \qquad \delta=2h,
\]
and compute the persistent-jump mass outside this band. For component $k$ we define
\[
    M^{\mathrm{tot}}_k
    =
    \sum_{x_j} P_{h,k}(x_j),
    \qquad
    M^{\mathrm{off}}_k
    =
    \sum_{x_j\notin B_\delta} P_{h,k}(x_j),
\]
summing over a diagnostic grid. The off-switch fraction is $\rho^{\mathrm{off}}_k=M^{\mathrm{off}}_k/M^{\mathrm{tot}}_k$, if $M^{\mathrm{tot}}_k>0$. 

\vspace{0.1in}

\paragraph{Rollout and vector field metrics}
Across the examples, the tables report trajectory-level and vector field-level
diagnostics.  Given a fixed set of $N$ initial conditions, the
rollout error is
\[
\mathrm{RMSE}_{\mathrm{rollout}}
=
\sqrt{
\frac{1}{N(K+1)}
\sum_{n=1}^{N}
\sum_{k=0}^{K}
\left\|
\hat x_n(t_k)-x_n(t_k)
\right\|_2^2
},
\]
and the final-time error is
\[
\mathrm{RMSE}_{\mathrm{final}}
=
\sqrt{
\frac{1}{N}
\sum_{n=1}^{N}
\left\|
\hat x_n(t_K)-x_n(t_K)
\right\|_2^2
}.
\]
These two quantities measure whether the learnt dynamics extrapolate beyond
the short training windows.

We also report vector field errors on a regular grid in the state domain,
averaged over forcing phases when the system is non-autonomous:
\[
\mathrm{RMSE}_{f}
=
\sqrt{
\frac{1}{M}
\sum_{m=1}^{M}
\left\|
\hat f(x_m,t_m)-f(x_m,t_m)
\right\|_2^2
}.
\]
The same quantity is also evaluated in a narrow band around the switching set.
This banded error is a strict local diagnostic of how well the switching dynamics are captured.

The evaluation protocol uses $128$ trajectories of $200$ steps for each trajectory-learning case.  Non-autonomous evaluation trajectories use independently sampled initial forcing phases.  Vector-field diagnostics use an $80\times 80$ state grid for dry friction and a $24^3$ grid for PP04. Each non-autonomous grid is averaged over $16$ uniformly spaced forcing phases. The switching band used for the local vector-field metric has half-width $0.05$.  Rollout metrics for the discontinuous dry-friction and PP04 models use the
single-hyperplane Filippov solver described in
Appendix~\ref{app:filippovSolver}; continuous Neural ODE baselines use RK4.

\vspace{0.1in}

\paragraph{Quality of the recovered switching geometry} When the switching geometry is learned, the tables report the angle and offset
errors of the recovered hyperplane.  For a learned hyperplane $a^\top x+b=0$
and the true hyperplane $a_\star^\top x+b_\star=0$, both
coefficient pairs are post-processed before evaluation.  The normals are rescaled to unit Euclidean norm and oriented so that their inner product is non-negative. The offsets are rescaled and, when necessary, sign-flipped with their corresponding normals.  Denote the resulting representations by $(\widehat a,\widehat b)$ and $(\widehat a_\star,\widehat b_\star)$.  The
reported angle error, in degrees, is
\[
\theta =
\frac{180}{\pi}
\arccos\left(
\widehat a^\top \widehat a_\star
\right),
\]
and the offset error is
\[
e_b =
\left|\widehat b-\widehat b_\star\right|.
\]
Both errors are therefore independent of the scaling and sign used to represent either hyperplane. For known-geometry runs, these quantities are zero by construction.

\subsection{Trajectory data generation}
\label{sec:trajectory_data_generation}
For each discontinuous-model trajectory experiment, we generate two independent sets of trajectories.  The main set is used for dynamics learning and is split into $90\%$ training and $10\%$ validation trajectories. Each dry-friction run uses a main set of $1,000$ trajectories, while each PP04 run uses $5,000$.  The warm-up set contains
the same number of trajectories as the corresponding main set and is used only to recover the switching geometry when that geometry is not supplied.  The continuous baselines use the same main trajectory sets and splits, but do not generate a warm-up set.

We sample uniformly the initial states from $[-3,3]\times[-2,2]$ for dry friction and from $[0,1]^3$ for PP04.  The initial time $t_0$ is sampled uniformly over one forcing period for each trajectory.

The sampled trajectory segment is
\[
    \{x(t_0),x(t_0+h),\ldots,x(t_0+kh)\}.
\]
For forced dry friction, reference trajectories are generated by a dedicated event-aware integrator.  It applies RK4 within each smooth branch, localises detected crossings of $p=0$ by bisection, and explicitly handles attractive sticking intervals on the switching manifold.  For PP04, the true flow is propagated exactly between events using the matrix exponential and a sinusoidal particular solution, see \cite[Section 6.2]{Morupisi2020}. Detected switching-plane crossings are localised by bisection before the remainder of the step is propagated on the new side.

Dry-friction main segments use $h=0.01$ and $k=10$, and warm-up segments use $h_{\rm warm}=0.01$ and $k_{\rm warm}=50$.  PP04 main segments use $h=0.05$ and $k=50$, and warm-up segments use $h_{\rm warm}=0.02$ and $k_{\rm warm}=100$.

\paragraph{Time supplied to the learned vector field}
The learned discontinuous networks do not receive raw time as an input.  Their smooth branch is evaluated as
\[
    \Psi\bigl(x,\Phi(z(x)),\tau(t)\bigr)
\]
when switching geometry is present, and as $\Psi(x,\tau(t))$ for a continuous model.  With a known forcing frequency, $\tau(t)=(\cos(\omega t),\sin(\omega t))$.  In the unknown-frequency experiment, $\tau(t)$ contains $16$ fixed random Fourier frequencies and phases, represented by their sine and cosine pairs. The universality result in \Cref{thm:approx-piecewise-smooth} persists within the family of piecewise-smooth vector fields, depending on $t$ only through $\tau(t)$.

\subsection{Known-geometry regression and discontinuity support}\label{sec:regression}

\begin{figure}[ht!]
    \centering
    \begin{subfigure}{0.96\linewidth}
        \centering
        \includegraphics[width=.8\linewidth]{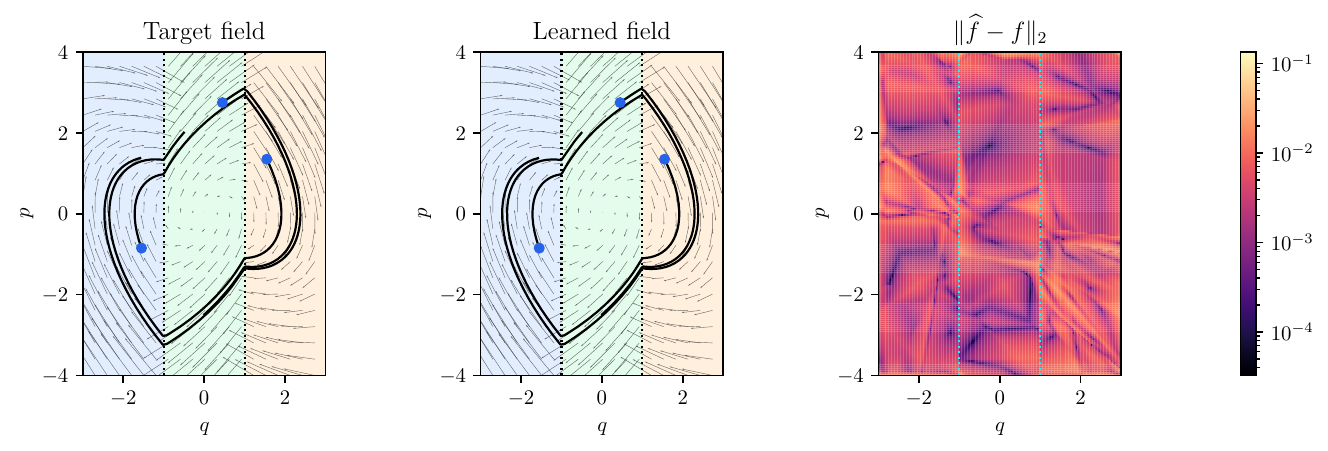}
        \caption{Geometry-preserving architecture.}
        \label{fig:regression_geometry_preserving_reconstruction}
    \end{subfigure}

    \vspace{0.75em}

    \begin{subfigure}{0.96\linewidth}
        \centering
        \includegraphics[width=.8\linewidth]{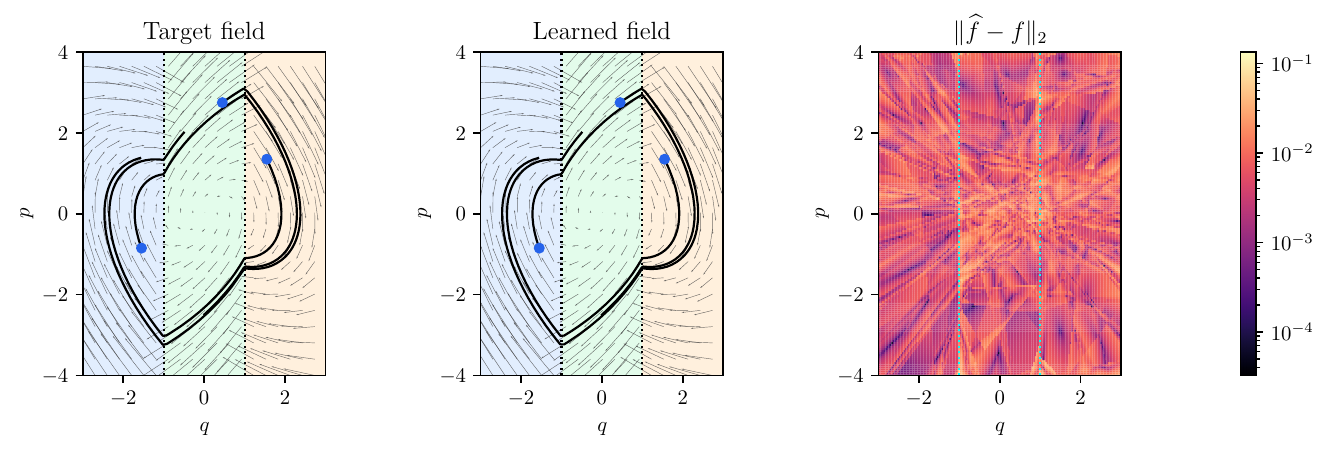}
        \caption{Unstructured jump network.}
        \label{fig:regression_unstructured_jump_reconstruction}
    \end{subfigure}
    \caption{Qualitative reconstruction of the target vector field from noisy
    regression data. In each row, the first panel shows the target field, the
    second panel shows the learned field, and the third panel shows the
    pointwise vector error $\|\widehat f-f\|_2$ on a logarithmic scale. Both
    architectures recover the main phase portrait and give small clean-grid
    errors, but the geometry-preserving architecture gives a smoother and more
    faithful reconstruction of the target partition structure. The
    displayed state domain is the full training rectangle
    $[-3,3]\times[-4,4]$, which contains every point of the diagnostic
    trajectories.}
    \label{fig:regression_reconstruction_comparison}
\end{figure}

We isolate the role of the architectural prior in the dynamics-learning problem by considering a two-dimensional piecewise affine vector field for which the switching geometry is known exactly, see \Cref{eq:pw-oscillator}. Both models are given the true switching code, and we ask whether the learned vector field introduces additional output jumps away from the prescribed switching set.

The switching set is given by:
\[
    \Sigma = \{(q,p)\in\R^2:\ q=-1\}\cup
    \{(q,p)\in\R^2:\ q=1\},
\]
which separates the state space into the three open polyhedra
\[
    \p_1=\{q<-1\},\qquad
    \p_2=\{-1<q<1\},\qquad
    \p_3=\{q>1\}.
\]
Equivalently, the implementation uses the binary switching code
\[
    z_P(x)=\mathcal H
    \left(
    \begin{bmatrix}
        1&0\\
        1&0
    \end{bmatrix}x+
    \begin{bmatrix}
        1\\
        -1
    \end{bmatrix}
    \right),
\]
so that the three regions correspond to the codes $(0,0)$, $(1,0)$, and $(1,1)$, respectively. 
\paragraph{Architectures}
The first model is the time-independent version of the geometry-preserving
class in \Cref{def:FP},
\[
    \mathcal N(x)=\Psi\bigl(x,\Phi(z_P(x))\bigr),
    \qquad
    \Phi\in\mathcal F_{\mathrm{ReLU}_J}^{2,d_1},
    \qquad
    \Psi\in\mathcal F_{\mathrm{ReLU}}^{2+d_1,2}.
\]
The two switching hyperplanes are fixed to their true values. We refer to this model as the geometry-preserving architecture.

We also compare against an unstructured jump network. This model keeps the same fixed switching code, but concatenates it directly with the state and passes the result through a jump-ReLU network,
\[
    \widetilde{\mathcal N}(x)
    =
    G\bigl(x,z_P(x)\bigr),
    \qquad
    G\in\mathcal F_{\mathrm{ReLU}_J}^{4,2}.
\]
This baseline therefore has access to the same exact partition information, but it does not use the separated $\Phi$--$\Psi$ parametrisation of \Cref{def:FP}. The widths are chosen so that the two models have comparable size.
\vspace{0.1in}

\paragraph{Training protocol}
Training data are generated directly from the vector field rather than from finite-difference trajectory estimates
\[
    x_n\sim\mathcal{U}\left([-3,3]\times[-4,4]\right),\,\,y_n = F(x_n) + \eta_n,\qquad
    \eta_n\sim\mathcal \mathsf{N}(0,\sigma^2 I_2),
    \qquad \sigma=0.03,
\]
for $n=1,\dots,N$. The model parameters are fitted by minimising the empirical mean-square regression loss
\[
    \mathcal{L}(\mathcal N):=\frac{1}{2N}\sum_{n=1}^N
    \left\|\mathcal N(x_n)-y_n\right\|_2^2.
\]

Each model is trained for $1,000$ epochs on $10,000$ noisy uniformly sampled points and selected using an independent noisy validation set of $2,000$ uniformly sampled points.  A separate clean test set contains $6,000$ uniformly sampled points.  The geometry-preserving model uses one width-$32$ hidden layer in $\Phi$ and three width-$32$ hidden layers in $\Psi$. The unstructured model uses three hidden layers of width $45$. The widths are chosen to give comparable model sizes: the geometry-preserving and unstructured models contain 4,482 and 4,592 trainable parameters, respectively. Both are optimised with Adam, weight decay $10^{-6}$, gradient-norm clipping at $1$, and a one-cycle learning-rate schedule from $10^{-4}$ to $5\cdot10^{-3}$ and back down by a final division factor of $10^4$.  The checkpoint with the smallest validation mean-square error is restored for all reported diagnostics.

For the orbit diagnostic, both the target and learned fields are propagated with classical RK4 on the same uniform output grid.  An outer RK4 step is recursively bisected into consecutive half-steps whenever the discrete region labels at its proposed endpoints differ.

\begin{figure}[ht!]
    \centering
    \includegraphics[width=0.6\linewidth]{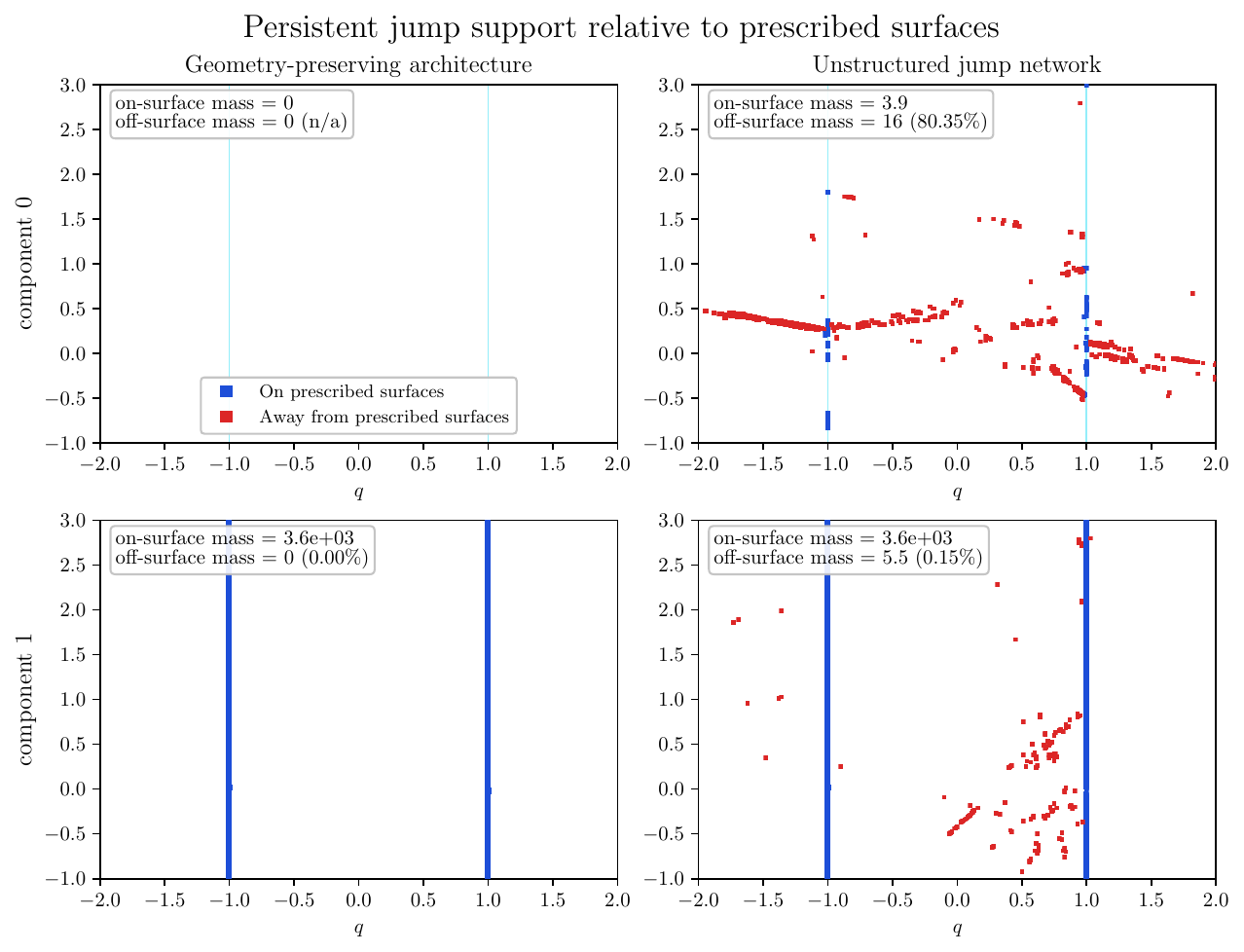}
    \caption{Persistent jump support after training on noisy regression targets
    with $\sigma=0.03$. Blue points lie within a $2h$ neighbourhood of the
    prescribed surfaces, and red points lie outside it. The architecture fixes
    the potential discontinuity surfaces but does not prescribe which output
    components may jump across them. No persistent jump mass is
    detected in the geometry-preserving model's first component, while all
    detected second-component mass lies on the prescribed surfaces.  The
    unstructured jump network produces persistent support away from those
    surfaces. Pale cyan lines in the first row mark potential architectural
    surfaces. The stronger solid guides in the second row mark the target
    component's actual switching surfaces.}
    \label{fig:regression_off_switch_jump_support}
\end{figure}

\paragraph{Results}
As shown in \Cref{fig:regression_reconstruction_comparison}, both models recover the main phase portrait from noisy observations.  The geometry-preserving model gives a more faithful field reconstruction and more closely reproduces the diagnostic trajectories throughout the training domain.  The persistent-jump diagnostic in \Cref{fig:regression_off_switch_jump_support} shows the central architectural difference. All detected persistent jumps produced by the geometry-preserving model remain confined to the prescribed switching surfaces, whereas the unstructured network also produces persistent jumps away from them.  Thus, access to the correct partition alone is not sufficient to preserve its geometry. The constrained architecture enforces the intended support for discontinuities while retaining an accurate approximation of the dynamics.

\subsection{Trajectory-based recovery of the discontinuities}\label{sec:learn-geometry}

\begin{figure}[ht!]
    \centering
    \includegraphics[width=\linewidth]{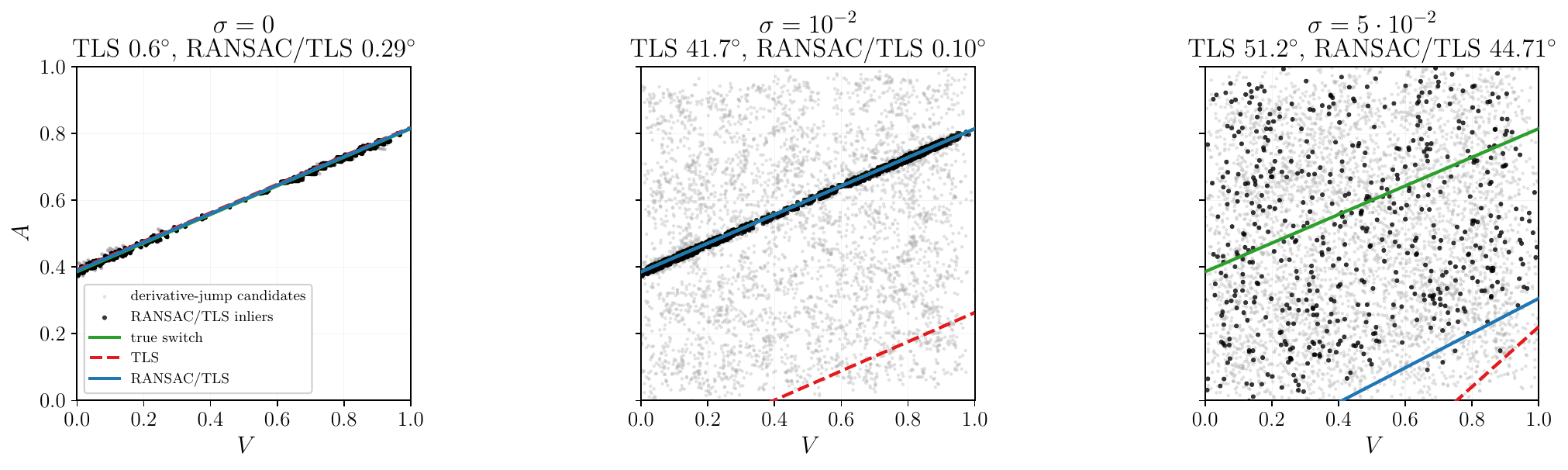}
    \caption{PP04 geometry recovery from derivative-jump candidates.  Grey
    points are the overcomplete candidate cloud formed from the endpoints of the
    globally retained high-derivative-jump intervals after neighbour
    expansion, and black points are the final RANSAC/TLS inliers.  Ordinary
    TLS and RANSAC receive the same candidate cloud.  Green is the
    true switch, red dashed is ordinary TLS, and blue is RANSAC followed by a
    final TLS refinement on its consensus set.  The panels correspond to
    observation-noise levels $\sigma=0$, $10^{-2}$, and $5\times10^{-2}$.
    The clean case shows that ordinary TLS works when the candidate cloud is
    already concentrated near the switch, the low-noise case shows the benefit
    of robust consensus selection in a contaminated cloud, and the high-noise
    case shows the recovery failure regime. At $\sigma=0$, the TLS
    and RANSAC angle errors are approximately $0.60^\circ$ and $0.29^\circ$.
    At $\sigma=10^{-2}$, TLS fails with error $41.7^\circ$, whereas RANSAC
    recovers the plane with error $0.10^\circ$.  At $\sigma=5\cdot10^{-2}$,
    both methods fail, with errors $51.2^\circ$ and $44.7^\circ$.}
    \label{fig:geometry-pp04-tls-vs-ransac}
\end{figure}

This subsection focuses on the geometry learning phase.  We consider the hyperplane learning procedure described in \Cref{alg:single_hyperplane_recovery_experiments,alg:ransac_experiments}. \Cref{fig:geometry-pp04-tls-vs-ransac} confirms that the candidate cloud of high-jump points that we work with after derivative estimation is overcomplete.  Large derivative jumps may come from true discontinuity crossings, but also from smooth high-curvature motion, endpoint effects in the derivative estimate, or observation noise. Consequently, fitting a hyperplane by total least squares (TLS) to all candidates is unreliable.  Instead, we use RANSAC to find a hyperplane supported by a coherent subset of candidates, and then refit it to the inliers using TLS.

Because both estimators receive exactly the same complete three-dimensional
candidate array, the comparison isolates the effect of robust consensus selection.  RANSAC selects the best-supported raw affine proposal, using the mean inlier distance only to break ties in support, and then applies TLS to the winning consensus.  This remains appropriate only when a coherent switching trace is identifiable in the candidate cloud. The comparison shows that robust consensus selection is decisive at mild noise but cannot recover an identifiable plane once the switching trace is not sufficiently distinguishable from the effects of noise.

\begin{figure}[ht!]
    \centering
    \includegraphics[width=0.5\linewidth]{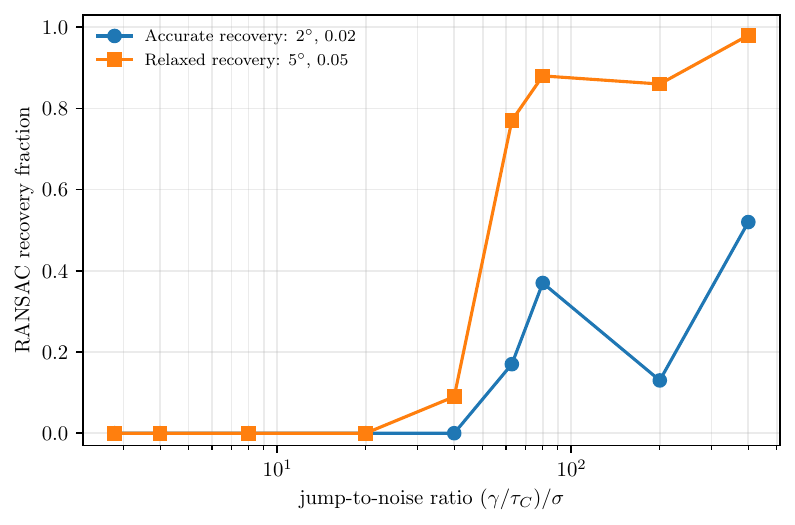}
    \caption{Empirical PP04 switching-plane recovery fraction as the
    discontinuity coefficient $\gamma$ is varied at fixed observation-noise
    level $\sigma=0.05$.  For each value of $\gamma$, the experiment generates
    $5,000$ trajectories and evaluates $100$ RANSAC seeds with $20,000$ proposals
    per seed.  Initial conditions, forcing phases, and standardised observation
    noise are paired across the values of $\gamma$.  The two curves report the
    fraction of seeds satisfying the accurate criterion (angle below $2^\circ$
    and offset error below $0.02$) and the relaxed criterion (angle below
    $5^\circ$ and offset error below $0.05$), respectively, against
    $(\gamma/\tau_C)/\sigma$.}
    \label{fig:geometry-pp04-jump-to-noise}
\end{figure}

We also empirically validate the detectability result in \Cref{sec:theory_geometryLearning}. We consider the PP04 geometry-recovery problem, varying the coefficient $\gamma$ which controls the magnitude of the jump across the discontinuity hyperplane. The observation-noise level remains $\sigma=0.05$.  Each success fraction is taken over 100 RANSAC random seeds for one paired trajectory/noise realisation. As shown in \Cref{fig:geometry-pp04-jump-to-noise}, recovery becomes
reliable once the discontinuity is sufficiently pronounced relative to the
observation noise, although the transition is not strictly monotone. We remark that changing
$\gamma$ also changes the trajectories and the resulting candidate cloud. Therefore, this test supports the theoretical analysis, but considers a different dynamical system for each $\gamma$.

\subsection{Trajectory-based learning of the geometry and the dynamics}\label{sec:full_pipeline}

The following experiments evaluate the problem of learning the dynamics from trajectory data after selecting a switching geometry.  In some runs, the true switching set is supplied to the model. In the remaining runs, it is recovered from warm-up trajectories and then kept fixed during training.  This separation is useful because geometry errors and vector field errors affect the learnt flow in different ways. When we use the trajectory rollout-loss to train our piecewise-smooth neural networks, we weigh the contributions to the MSE loss coming from the two sides of the switch so they are balanced. The balancing is based on the sign of the switching function at the trajectory initial condition.

\subsubsection{Dry-friction oscillators}

The dry-friction examples, see \Cref{eq:dry-friction}, provide a controlled setting in which the true
switching set is simple, but the long-time behaviour is sensitive to both the
discontinuity and the sliding dynamics.  The observed trajectories are perturbed with additive Gaussian noise with $\sigma\in\{0,10^{-2},5\cdot 10^{-2}\}$.
When the second-order structure is enforced, we only parametrise with the network the second component of the vector field, and force the first to be equal to $p$. For each case and noise level, we train the four combinations assuming the knowledge, or lack thereof, of the switching hyperplane and second-order nature of the dynamics. The forcing frequency is assumed known in the main dry-friction tables. Appendix \ref{app:additionalDryp} shows that replacing the known harmonic features with fixed random Fourier features substantially degrades both local accuracy and temporal extrapolation. All dry-friction models use one hidden layer of width $10$ in $\Phi$ and
three hidden layers of width $10$ in $\Psi$. In the learned-switch runs, we
retain the largest $2\%$ of the derivative-jump scores globally, augment each
retained interval with its neighbours within radius $1$, and use both interval
endpoints as candidates. We fit one switching line using $20,000$ RANSAC
proposals with distance tolerance $0.03$, followed by total least-squares
refinement on the inliers. All models minimise \eqref{eq:lossFull}.  The velocities
are estimated with a Savitzky--Golay window of $11$ and polynomial degree $3$.
The loss weights are $\lambda_{\rm traj}=1$,
$\lambda_{\rm vf}=10^{-2}$, $\lambda_{\rm nor}=3\cdot10^{-3}$, and
$\lambda_{\rm tan}=10^{-2}$.  The velocity loss excludes
$|s(x)|<0.05$, and the normal and tangential thresholds are
$\tau_{\rm nor}=5\cdot10^{-3}$ and $\tau_{\rm tan}=10^{-3}$.

Training uses classical RK4, not the Filippov solver. This choice is based on empirical evidence of comparable performance but at a much lower time cost. The rollout horizon follows a curriculum from $2$ steps to the full $10$-step window.  Adam is run for $500$ epochs with batch size $128$, weight decay
$10^{-5}$, gradient-norm clipping at $1$, and a one-cycle schedule between
$5\cdot10^{-4}$ and $5\cdot10^{-3}$.  Full-horizon training and validation
losses are measured every five epochs and after the final epoch, and the
checkpoint with the smallest validation rollout loss is restored.  The Filippov solver, see Appendix~\ref{app:filippovSolver}, is used for rollout evaluation of the discontinuous models.

The recovered line is fixed during dynamics learning.  In the reported tables, the column ``2nd order''
indicates whether the model is constrained to respect $\dot q=p$, so that only
the acceleration component is learnt.

\begin{table*}[ht!]
\centering
\caption{Dynamics-learning results for dry-friction case 1. Each sub-table corresponds to a different trajectory-noise level.}
\label{tab:traj-dryp_case1_learning_dynamics}
\begin{subtable}{\textwidth}
\centering
\caption{Clean data, $\sigma=0$}
\label{tab:traj-dryp_case1_noise0}
\scriptsize
\resizebox{\textwidth}{!}{%
\begin{tabular}{llrrrrrr}
\toprule
Known switch & 2nd order & Switch angle & Switch offset & Rollout RMSE & Final RMSE & VF RMSE & VF band RMSE \\
\midrule
no & no & 0.225 & $3.78\cdot 10^{-3}$ & 0.0348 & 0.0701 & 0.0146 & $3.00\cdot 10^{-3}$ \\
no & yes & 0.225 & $3.78\cdot 10^{-3}$ & $3.90\cdot 10^{-3}$ & $5.17\cdot 10^{-3}$ & $6.27\cdot 10^{-3}$ & $4.65\cdot 10^{-3}$ \\
yes & no & 0 & 0 & 0.0356 & 0.0719 & 0.0144 & $2.86\cdot 10^{-3}$ \\
yes & yes & 0 & 0 & $3.18\cdot 10^{-3}$ & $3.99\cdot 10^{-3}$ & $6.41\cdot 10^{-3}$ & $4.73\cdot 10^{-3}$ \\
\bottomrule
\end{tabular}%
}
\end{subtable}

\vspace{0.6em}
\begin{subtable}{\textwidth}
\centering
\caption{Mild noise, $\sigma=10^{-2}$}
\label{tab:traj-dryp_case1_noise001}
\scriptsize
\resizebox{\textwidth}{!}{%
\begin{tabular}{llrrrrrr}
\toprule
Known switch & 2nd order & Switch angle & Switch offset & Rollout RMSE & Final RMSE & VF RMSE & VF band RMSE \\
\midrule
no & no & 0.284 & $4.92\cdot 10^{-3}$ & 0.0882 & 0.156 & 0.0618 & 0.0613 \\
no & yes & 0.284 & $4.92\cdot 10^{-3}$ & 0.0150 & 0.0235 & 0.0202 & 0.0231 \\
yes & no & 0 & 0 & 0.0887 & 0.157 & 0.0623 & 0.0620 \\
yes & yes & 0 & 0 & 0.0219 & 0.0351 & 0.0291 & 0.0315 \\
\bottomrule
\end{tabular}%
}
\end{subtable}

\vspace{0.6em}
\begin{subtable}{\textwidth}
\centering
\caption{High noise, $\sigma=5\cdot 10^{-2}$}
\label{tab:traj-dryp_case1_noise005}
\scriptsize
\resizebox{\textwidth}{!}{%
\begin{tabular}{llrrrrrr}
\toprule
Known switch & 2nd order & Switch angle & Switch offset & Rollout RMSE & Final RMSE & VF RMSE & VF band RMSE \\
\midrule
no & no & 1.34 & 0.0475 & 0.202 & 0.331 & 0.227 & 0.553 \\
no & yes & 1.34 & 0.0475 & 0.0869 & 0.129 & 0.164 & 0.514 \\
yes & no & 0 & 0 & 0.275 & 0.446 & 0.259 & 0.276 \\
yes & yes & 0 & 0 & 0.125 & 0.184 & 0.193 & 0.220 \\
\bottomrule
\end{tabular}%
}
\end{subtable}
\end{table*}

\begin{table*}[ht!]
\centering
\caption{Dynamics-learning results for dry-friction case 2. Each sub-table corresponds to a different trajectory-noise level.}
\label{tab:traj-dryp_case2_learning_dynamics}
\begin{subtable}{\textwidth}
\centering
\caption{Clean data, $\sigma=0$}
\label{tab:traj-dryp_case2_noise0}
\scriptsize
\resizebox{\textwidth}{!}{%
\begin{tabular}{llrrrrrr}
\toprule
Known switch & 2nd order & Switch angle & Switch offset & Rollout RMSE & Final RMSE & VF RMSE & VF band RMSE \\
\midrule
no & no & 0.194 & $9.58\cdot 10^{-4}$ & 0.0234 & 0.0537 & 0.0146 & $3.18\cdot 10^{-3}$ \\
no & yes & 0.194 & $9.58\cdot 10^{-4}$ & $2.16\cdot 10^{-3}$ & $3.16\cdot 10^{-3}$ & $6.25\cdot 10^{-3}$ & $4.70\cdot 10^{-3}$ \\
yes & no & 0 & 0 & 0.0236 & 0.0543 & 0.0144 & $2.93\cdot 10^{-3}$ \\
yes & yes & 0 & 0 & $2.17\cdot 10^{-3}$ & $3.40\cdot 10^{-3}$ & $5.89\cdot 10^{-3}$ & $4.61\cdot 10^{-3}$ \\
\bottomrule
\end{tabular}%
}
\end{subtable}

\vspace{0.6em}
\begin{subtable}{\textwidth}
\centering
\caption{Mild noise, $\sigma=10^{-2}$}
\label{tab:traj-dryp_case2_noise001}
\scriptsize
\resizebox{\textwidth}{!}{%
\begin{tabular}{llrrrrrr}
\toprule
Known switch & 2nd order & Switch angle & Switch offset & Rollout RMSE & Final RMSE & VF RMSE & VF band RMSE \\
\midrule
no & no & 0.300 & $2.54\cdot 10^{-3}$ & 0.0748 & 0.138 & 0.0617 & 0.0613 \\
no & yes & 0.300 & $2.54\cdot 10^{-3}$ & 0.0214 & 0.0343 & 0.0286 & 0.0311 \\
yes & no & 0 & 0 & 0.0752 & 0.139 & 0.0622 & 0.0618 \\
yes & yes & 0 & 0 & 0.0224 & 0.0356 & 0.0294 & 0.0318 \\
\bottomrule
\end{tabular}%
}
\end{subtable}

\vspace{0.6em}
\begin{subtable}{\textwidth}
\centering
\caption{High noise, $\sigma=5\cdot 10^{-2}$}
\label{tab:traj-dryp_case2_noise005}
\scriptsize
\resizebox{\textwidth}{!}{%
\begin{tabular}{llrrrrrr}
\toprule
Known switch & 2nd order & Switch angle & Switch offset & Rollout RMSE & Final RMSE & VF RMSE & VF band RMSE \\
\midrule
no & no & 1.29 & 0.0332 & 0.170 & 0.293 & 0.191 & 0.514 \\
no & yes & 1.29 & 0.0332 & 0.131 & 0.196 & 0.207 & 0.484 \\
yes & no & 0 & 0 & 0.200 & 0.335 & 0.195 & 0.180 \\
yes & yes & 0 & 0 & 0.126 & 0.186 & 0.195 & 0.222 \\
\bottomrule
\end{tabular}%
}
\end{subtable}
\end{table*}
\paragraph{Results}
\Cref{tab:traj-dryp_case1_learning_dynamics} and \Cref{tab:traj-dryp_case2_learning_dynamics} show that the second-order constraint improves
rollout and final-state RMSE in every case, noise level, and geometry setting.
At $\sigma\leq10^{-2}$, the recovered switches remain accurate, with angle
errors below $0.30^\circ$ and offset errors below $5\cdot10^{-3}$.  At
$\sigma=5\cdot10^{-2}$, the angle errors increase to approximately
$1.3^\circ$ and the switching-band errors increase substantially.  Known
geometry does not uniformly improve rollout error. The learned-switch and
known-switch models are optimised independently, and the recovered geometry
can act as a slightly different inductive bias.  The consistent result is the
benefit of second-order structure and the degradation of local switching-band
accuracy when high-noise geometry recovery deteriorates.

\subsubsection{PP04}\label{sec:traj-pp04_experiment}

We use PP04, see \Cref{eq:pp04}, as
a second trajectory-learning benchmark after the dry-friction examples.  The
switching geometry must either be supplied to the model or recovered from
trajectory data, and the learned vector field is then assessed through both
local vector field errors and finite-horizon rollouts.

The forcing frequency $\omega$ is assumed known in all PP04 runs. 
The observed trajectories have initial conditions sampled in $[0,1]^3$ and perturbed by additive Gaussian noise with standard deviation $\sigma\in\{0,10^{-2}\}$. We choose these noise levels because of the analysis in \Cref{sec:theory_geometryLearning}. All PP04 models use one hidden layer of width $11$ in $\Phi$ and one hidden
layer of width $14$ in $\Psi$. In the learned-switch runs, we retain the largest
$0.5\%$ of derivative-jump scores globally, augment each retained interval with
its neighbours within radius $2$, and use both interval endpoints as candidates.
We fit one switching plane using $20,000$ RANSAC proposals with distance
tolerance $10^{-2}$, followed by total least-squares refinement on the inliers.

Dynamics training minimises \eqref{eq:lossFull} using classical RK4,
with $\lambda_{\rm traj}=1$, $\lambda_{\rm vf}=10^{-2}$, no switch-band
exclusion, and $\lambda_{\rm nor}=\lambda_{\rm tan}=0$. Dynamics-loss velocities use a Savitzky--Golay window of
$21$ and polynomial degree $3$; warm-up geometry recovery instead uses window
$51$ and degree $2$.  The rollout curriculum grows from $2$ to $50$ steps.
Adam is run for $500$ epochs with batch size $128$, weight decay $10^{-5}$,
gradient-norm clipping at $1$, and a one-cycle schedule between $10^{-5}$ and
$10^{-4}$.  Validation is evaluated every five epochs and after the final epoch,
and the best full-horizon validation checkpoint is restored.
\begin{table*}[ht!]
\centering
\caption{PP04 dynamics-learning results with known forcing frequency and Filippov rollout evaluation. Each sub-table fixes the noise level and compares the full learned-geometry pipeline against the known-switch setting.}
\label{tab:traj-pp04_learning_dynamics}
\begin{subtable}{\linewidth}
\centering
\caption{$\sigma=0$}
\label{tab:traj-pp04_learning_noise_0p0}
\small
\begin{tabular}{lrrrrrr}
\toprule
\makecell{Known\\switch} & \makecell{Switch\\angle} & \makecell{Switch\\offset} & \makecell{Rollout\\RMSE} & \makecell{Final\\RMSE} & \makecell{VF\\RMSE} & \makecell{VF band\\RMSE} \\
\midrule
no & 0.339 & $5.18\cdot 10^{-3}$ & $2.66\cdot 10^{-3}$ & $2.67\cdot 10^{-3}$ & $6.18\cdot 10^{-3}$ & 0.0191 \\
yes & 0 & 0 & $1.45\cdot 10^{-3}$ & $1.7\cdot 10^{-3}$ & $1.3\cdot 10^{-3}$ & $2.19\cdot 10^{-3}$ \\
\bottomrule
\end{tabular}
\end{subtable}
\vspace{0.75em}
\begin{subtable}{\linewidth}
\centering
\caption{$\sigma=10^{-2}$}
\label{tab:traj-pp04_learning_noise_0p01}
\small
\begin{tabular}{lrrrrrr}
\toprule
\makecell{Known\\switch} & \makecell{Switch\\angle} & \makecell{Switch\\offset} & \makecell{Rollout\\RMSE} & \makecell{Final\\RMSE} & \makecell{VF\\RMSE} & \makecell{VF band\\RMSE} \\
\midrule
no & 0.458 & $4.37\cdot 10^{-3}$ & 0.0111 & 0.0176 & $7.9\cdot 10^{-3}$ & 0.0244 \\
yes & 0 & 0 & $8.07\cdot 10^{-3}$ & 0.0125 & $1.68\cdot 10^{-3}$ & $2\cdot 10^{-3}$ \\
\bottomrule
\end{tabular}
\end{subtable}
\end{table*}

The recovered switching plane remains accurate at both noise
levels.  All four models reproduce the PP04 dynamics closely, although
observation noise increases both rollout and local vector-field errors.
Supplying the true switching plane reduces every reported error at each noise
level, with the largest improvement occurring near the switching surface.

The initial state is shared across panels and selected from $30$
deterministic candidates by minimising aggregate rollout error over
$0\leq t\leq 250$.  Every candidate rollout starts at $t_0=0$.  The figure displays the first $500$ time units
for readability.  This selection is used only for visualisation and does not
affect training, checkpoint selection, or the evaluation metrics.

\begin{figure}[ht!]
    \centering
    \includegraphics[width=\linewidth]{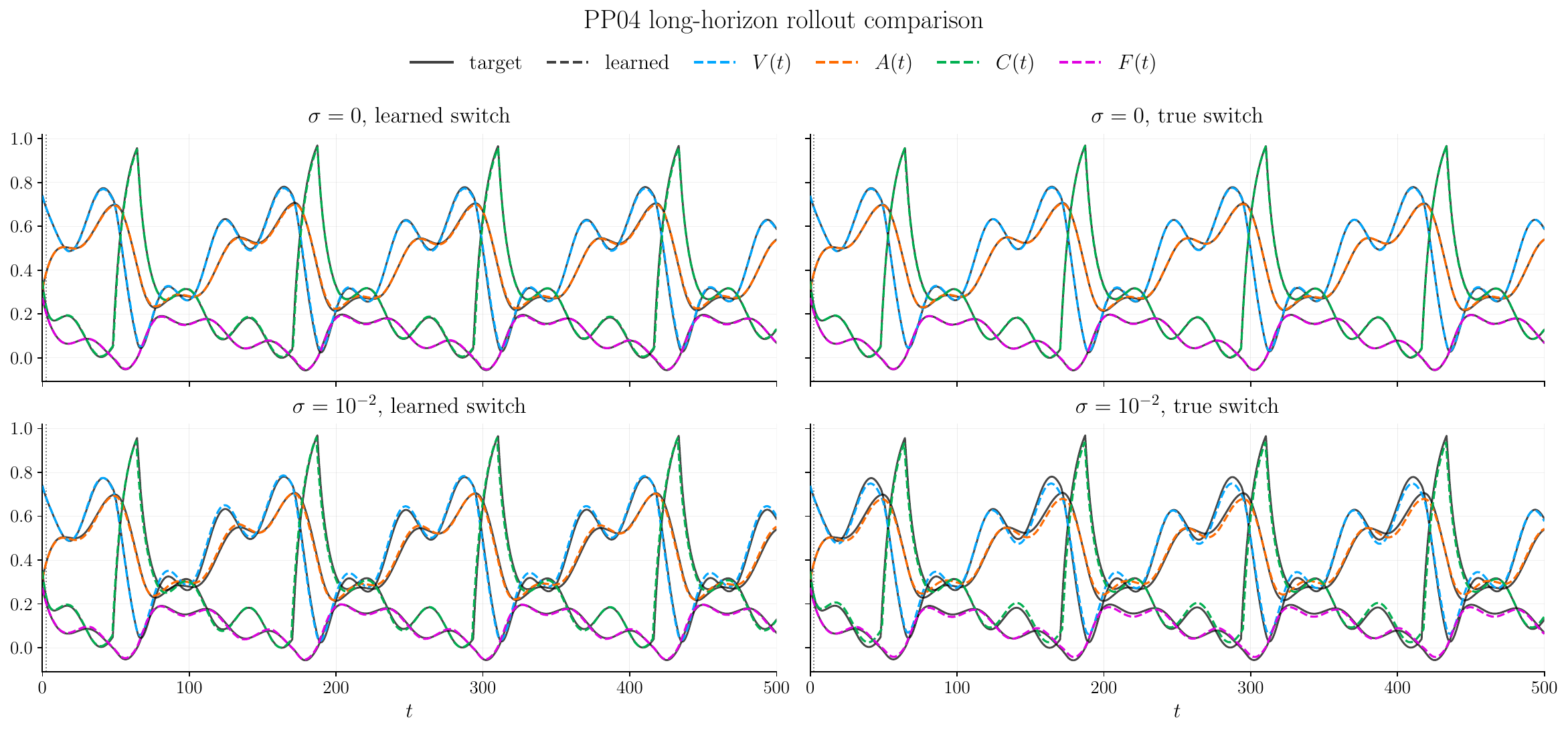}
    \caption{PP04 rollouts over $0\leq t\leq 500$ from the shared
    initial state selected by minimum aggregate error across all four models,
    with $t_0=0$. Each panel corresponds to one noise and switching-geometry
    setting, and each trace overlays the target trajectory with the learnt
    trajectory for $V$, $A$, $C$, and the switching expression
    $aV-bA+d$.  The PP04 state itself has three components. The clean runs are nearly indistinguishable from the target,
    while the noisy runs retain the recurrent structure with small local
    deviations.}
    \label{fig:traj-pp04-rollout-grid}
\end{figure}

\subsubsection{Comparison with a smooth Neural ODE}
We compare the structured models with a smooth Neural ODE trained on
the same trajectory pools, splits, segment lengths, and known-frequency time
features.  The baseline is a three-hidden-layer width-$64$ MLP with $\tanh$
activations.  It is trained for $500$ epochs using Adam at the fixed learning
rate $10^{-3}$, weight decay $10^{-6}$, batch size $128$ for dry friction and for PP04, and full-horizon RK4 trajectory mean-square error. Validation
is measured every five epochs and after the final epoch, and the best validation
checkpoint is restored.  The baseline has no switching code, geometry-recovery
stage, vector-field loss, or directional regulariser. This choice follows from empirically seeing that these additional components did not improve the model performance. The trainable parameter counts are constant across noise levels and
across the two dry-friction forcing cases.  For dry friction, the structured
\texttt{Disc} model has $521$ trainable parameters and the smooth
\texttt{NODE} baseline has $8705$.  For PP04, the corresponding counts are
$448$ and $8899$.  These totals count only parameters with gradients, not the fixed
switching-hyperplane tensors.
\begin{table}[ht!]
\centering
\caption{Comparison with a smooth Neural ODE baseline.  ``Disc'' denotes the learned-geometry discontinuous model and ``NODE'' the continuous baseline. Case 1 and 2 are dry-friction oscillators. For dry friction, both the discontinuous model and the smooth Neural ODE impose the second-order relation $\dot q=p$ and parametrise only $\dot p$. For PP04, both models parametrise all three components of the vector field. ``VF band'' restricts vector-field error to the switching band.  Bold indicates the smaller value within each matched pair. Discontinuous-model rollouts use the Filippov solver, while continuous
baselines use RK4.}
\label{tab:vanilla-node-comparison}
\begin{tabular}{llcccc}
\toprule
System & Model & \makecell{Rollout\\RMSE} &  \makecell{Final\\RMSE} & \makecell{VF\\RMSE} & \makecell{VF band\\RMSE} \\
\midrule
Case 1, $\sigma=0$
& Disc & $\mathbf{3.90\cdot 10^{-3}}$ & $\mathbf{5.17\cdot 10^{-3}}$ & $\mathbf{6.27\cdot 10^{-3}}$ & $\mathbf{4.65\cdot 10^{-3}}$ \\
Case 1, $\sigma=0$
& NODE & $3.53\cdot 10^{-2}$ & $5.55\cdot 10^{-2}$ & $8.65\cdot 10^{-2}$ & $4.17\cdot 10^{-1}$ \\
\addlinespace
Case 1, $\sigma=10^{-2}$
& Disc & $\mathbf{1.50\cdot 10^{-2}}$ & $\mathbf{2.35\cdot 10^{-2}}$ & $\mathbf{2.02\cdot 10^{-2}}$ & $\mathbf{2.31\cdot 10^{-2}}$ \\
Case 1, $\sigma=10^{-2}$
& NODE & $4.88\cdot 10^{-2}$ & $7.51\cdot 10^{-2}$ & $1.03\cdot 10^{-1}$ & $4.34\cdot 10^{-1}$ \\
\addlinespace
Case 2, $\sigma=0$
& Disc & $\mathbf{2.16\cdot 10^{-3}}$ & $\mathbf{3.16\cdot 10^{-3}}$ & $\mathbf{6.25\cdot 10^{-3}}$ & $\mathbf{4.70\cdot 10^{-3}}$ \\
Case 2, $\sigma=0$
& NODE & $2.79\cdot 10^{-2}$ & $4.43\cdot 10^{-2}$ & $8.65\cdot 10^{-2}$ & $4.17\cdot 10^{-1}$ \\
\addlinespace
Case 2, $\sigma=10^{-2}$
& Disc & $\mathbf{2.14\cdot 10^{-2}}$ & $\mathbf{3.43\cdot 10^{-2}}$ & $\mathbf{2.86\cdot 10^{-2}}$ & $\mathbf{3.11\cdot 10^{-2}}$ \\
Case 2, $\sigma=10^{-2}$
& NODE & $4.52\cdot 10^{-2}$ & $6.75\cdot 10^{-2}$ & $1.03\cdot 10^{-1}$ & $4.34\cdot 10^{-1}$ \\
\addlinespace
PP04, $\sigma=0$ & Disc & $\mathbf{2.66\cdot 10^{-3}}$ & $\mathbf{2.67\cdot 10^{-3}}$ & $\mathbf{6.18\cdot 10^{-3}}$ & $\mathbf{1.91\cdot 10^{-2}}$ \\
PP04, $\sigma=0$ & NODE & $3.28\cdot 10^{-2}$ & $6.24\cdot 10^{-2}$ & $8.96\cdot 10^{-3}$ & $2.71\cdot 10^{-2}$ \\
\addlinespace
PP04, $\sigma=10^{-2}$ & Disc & $\mathbf{1.11\cdot 10^{-2}}$ & $\mathbf{1.76\cdot 10^{-2}}$ & $\mathbf{7.9\cdot 10^{-3}}$ & $\mathbf{2.44\cdot 10^{-2}}$ \\
PP04, $\sigma=10^{-2}$ & NODE & $2.61\cdot 10^{-2}$ & $3.32\cdot 10^{-2}$ & $1.05\cdot 10^{-2}$ & $3.11\cdot 10^{-2}$ \\
\bottomrule
\end{tabular}
\end{table}

The compact discontinuous model outperforms the substantially larger NODE baseline in every metric at both reported noise levels. Thus the comparison supports strong parameter efficiency and model flexibility for the considered examples.

\section{Discussion and conclusion}
We have proposed a two-stage framework for learning piecewise-smooth dynamical systems from trajectory data. The method first estimates affine switching geometry from derivative jumps, and then learns dynamics adapted to the regions induced by the recovered or prescribed switching set. This separation makes the model interpretable and allows prior knowledge to be incorporated either in the geometry-learning stage or in the local dynamics.

The theoretical analysis provides results for both stages. It derives a local expected detectability condition that shows how observation noise, derivative-estimation error, smoothing, and jump size interact under regularity and transversal-crossing assumptions. It also shows that the proposed geometry-preserving neural class can uniformly approximate piecewise-smooth vector fields on compact sets away from the prescribed switching hyperplanes, while confining discontinuities to the prescribed partition.

The experiments support these findings. With known geometry, the proposed architecture avoids spurious off-switch jumps. With the recovered geometry, the learned models produce accurate rollouts when the derivative-jump signal is sufficiently strong relative to the noise. The dry-friction and PP04 tests also show the value of incorporating structural information, such as second-order dynamics or known forcing frequencies. More broadly, the results show that explicitly learning and preserving switching geometry offers a practical route to interpretable models of discontinuous dynamics.

The present study is restricted to fully observed low-dimensional states, affine switching geometry, and a sequential geometry--dynamics pipeline. The Filippov-aware solver is further restricted to a single switching hyperplane.

Future work will consider experimental noisy data, higher-dimensional systems, interacting switches, and non-affine switching surfaces. Another natural direction is to combine the method with encoder--decoder models, so that discontinuous latent dynamics can be learned from videos or other high-dimensional scientific observations.

\section*{Acknowledgements}
The authors would like to thank Dr Alice Cicirello for the insightful discussions at the start of this project. EJ was funded by the Knut and Alice Wallenberg Foundation grant 2024.0440. CBS acknowledges support from the Royal Society Wolfson Fellowship, the EPSRC advanced career fellowship EP/V029428/1, the Wellcome Innovator Awards 215733/Z/19/Z and 221633/Z/20/Z. DM, CB, and CBS acknowledge the support of the EPSRC programme grant EP/V026259/1. DM and CBS acknowledge support from the EPSRC programme grant EP/Y028783/1. All the authors acknowledge support from the EU through the Marie Sk\l{}odowska-Curie Actions Staff Exchanges project REMODEL, grant agreement 101131557.

 \paragraph{AI use declaration}
ChatGPT was used to assist in checking mathematical proofs and identifying notation issues and inconsistencies across the manuscript prior to submission. Codex was used to assist with the preparation of the research code, under direct and substantial supervision by the authors. The code is publicly available in the GitHub repository associated with this paper. The authors assume responsibility for all content.

\bibliographystyle{siamplain}
\bibliography{biblio}

\newpage

\appendix
\section{Fitting linear models in the presence of outliers} \label{app:ransaclinreg}

Out of a sample consisting of $M$ total points, with $M_i$ inliers and $M_o$ outliers, RANSAC is a randomised algorithm that iteratively separates inliers from outliers.  It works by randomly selecting $d$ points to determine a $(d-1)$-dimensional affine subspace of $\mathbb{R}^d$.
Then, a hyperplane is fitted to these points, and the points sufficiently close (say, within a distance $\varepsilon$) are included in the \emph{consensus set}. A model is deemed as reasonable if the consensus set is sufficiently large. The algorithm is repeated for a fixed number of iterations, each time keeping the model with the largest consensus set. 

The number of inliers selected follows a hypergeometric distribution (as there is no replacement), so the probability $r$ of selecting exactly $d$ inliers is given by 
\[
r = \frac{{\binom{M_i}{d}}}{{\binom{M}{d}}}.
\]
The probability of selecting at least one point which is not an inlier is therefore $1-r$, and the probability of this happening $k$ times is thus $(1-r)^k$, which is the probability that the RANSAC algorithm has not found a good model in $k$ iterations. The probability that at least one iteration has included only inliers is therefore $p = 1-(1-r)^k$. 
Thus, to ensure that RANSAC has reached the desired accuracy in $k$ iterations with probability $p$, we must take
\[
k \geq  \frac{\log(1-p)}{\log(1-r)}.
\] 
We see that we must estimate $r$, meaning we need to determine, or at least bound, the number of points included as candidates for the RANSAC algorithm, that is, $M$ and $M_i$. 

\Cref{alg:ransac_experiments} describes the hyperplane-fitting step in the geometry learning phase.

\begin{algorithm}[ht!]
\caption{$\operatorname{RANSAC}(\mathcal C,\varepsilon,M,m_{\rm req})$}
\label{alg:ransac_experiments}
\begin{algorithmic}[1]
\Require Candidate cloud $\mathcal C\subset\mathbb{R}^d$; tolerance $\varepsilon$; iterations $M$; required support $m_{\mathrm{req}}$.
\Ensure Best-supported hyperplane proposal and its consensus set, or failure.
\State Initialise the current winner as empty.
\For{$m=1,\dots,M$}
    \State Sample $d$ points uniformly without replacement from $\mathcal C$.
    \State Construct the hyperplane
    $\Gamma_m=\{x\in\mathbb R^d:a_m^\top x+b_m=0\}$ through the sampled points.
    \If{the sampled points are degenerate}
        \State \textbf{continue}.
    \EndIf
    \State Compute the consensus set $\mathcal I_m
    =
    \{x\in\mathcal C:\operatorname{dist}(x,\Gamma_m)\leq\varepsilon\}$.
    \If{$|\mathcal I_m|<m_{\rm req}$}
        \State \textbf{continue}.
    \EndIf
    \State Compute the mean inlier distance
    \[
    D_m=
    \frac{1}{|\mathcal I_m|}
    \sum_{x\in\mathcal I_m}\operatorname{dist}(x,\Gamma_m).
    \]
    \State Replace the winner if
    $(|\mathcal I_m|,-D_m)$ is lexicographically larger than its current score.
\EndFor
\If{no winner was found}
    \State \textbf{fail}.
\EndIf
\State \Return $(a_\star,b_\star,\mathcal I_\star)$.
\end{algorithmic}
\end{algorithm}

\section{Additional proofs for the universality result}\label{app:additionalProofsUniversal}

\begin{proof}[Proof of \Cref{lem:closure-basic}]
The arguments in this proof apply indistinguishably to the elements in $\mathcal F_{\mathrm{ReLU}}^{d_1,d_2}$ and $\mathcal F_{\mathrm{ReLU}_J}^{d_1,d_2}$. Scalar multiplication and post-composition by affine maps are immediate. For addition and parallel concatenation, one first pads the shallower network with identity layers of the form \eqref{eq:xJump-note} so that the two networks have the same depth. One may then use block-diagonal affine maps and componentwise activations to run the two networks in parallel. A final affine map gives either the concatenation or the sum.
\end{proof}

\begin{proof}[Proof of \Cref{prop:FP-linear}]
Let $f,g\in\mathcal F_P^d$. By definition there are dimensions $r_f,r_g\in\mathbb N$, jump networks $\Phi_f\in\mathcal F_{{\mathrm{ReLU}}_J}^{d_P,r_f}$ and $\Phi_g\in\mathcal F_{\mathrm{ReLU}_J}^{d_P,r_g}$, and continuous networks $\Psi_f\in\mathcal F_{\mathrm{ReLU}}^{d+r_f+1,d}$ and $\Psi_g\in\mathcal F_{\mathrm{ReLU}}^{d+r_g+1,d}$ such that
\[
f(x,t)=\Psi_f\bigl(x,t,\Phi_f(z_P(x))\bigr),
\qquad
g(x,t)=\Psi_g\bigl(x,t,\Phi_g(z_P(x))\bigr).
\]
By \Cref{lem:closure-basic}, the parallel concatenation $\Phi(z):=\bigl[\Phi_f(z),\Phi_g(z)\bigr]\in\mathcal F_{\mathrm{ReLU}_J}^{d_P,r_f+r_g}$ and the parallel concatenation of $\Psi_f$ and $\Psi_g$ both belong to the corresponding ReLU classes. A final affine map produces any linear combination of $f$ and $g$. Hence $\mathcal F_P^d$ is a vector space.
\end{proof}

\section{Additional experiments for the dry friction system}

\subsection{Representative dry-friction rollouts}
\label{app:dryp-representative-rollouts}

The aggregate errors in
\Cref{tab:traj-dryp_case1_learning_dynamics,tab:traj-dryp_case2_learning_dynamics}
summarise performance over the evaluation set.  Here we show representative
long-horizon trajectories to illustrate the qualitative behaviour behind those
numbers.  We fix the initial condition $(q_0,p_0)=(0.1,0.1)$ and $t_0=0$ for all runs. All panels use the second-order model, known forcing frequency,
noise level $\sigma=10^{-2}$, and the Filippov-aware rollout solver.  This
noise level is useful for visualisation because it is large enough to make
geometry recovery nontrivial, but not so large that the recovered-switch model
is dominated by geometry failure.

\begin{figure}[ht!]
    \centering
    \begin{subfigure}{0.49\linewidth}
        \centering
        \includegraphics[width=.9\linewidth]{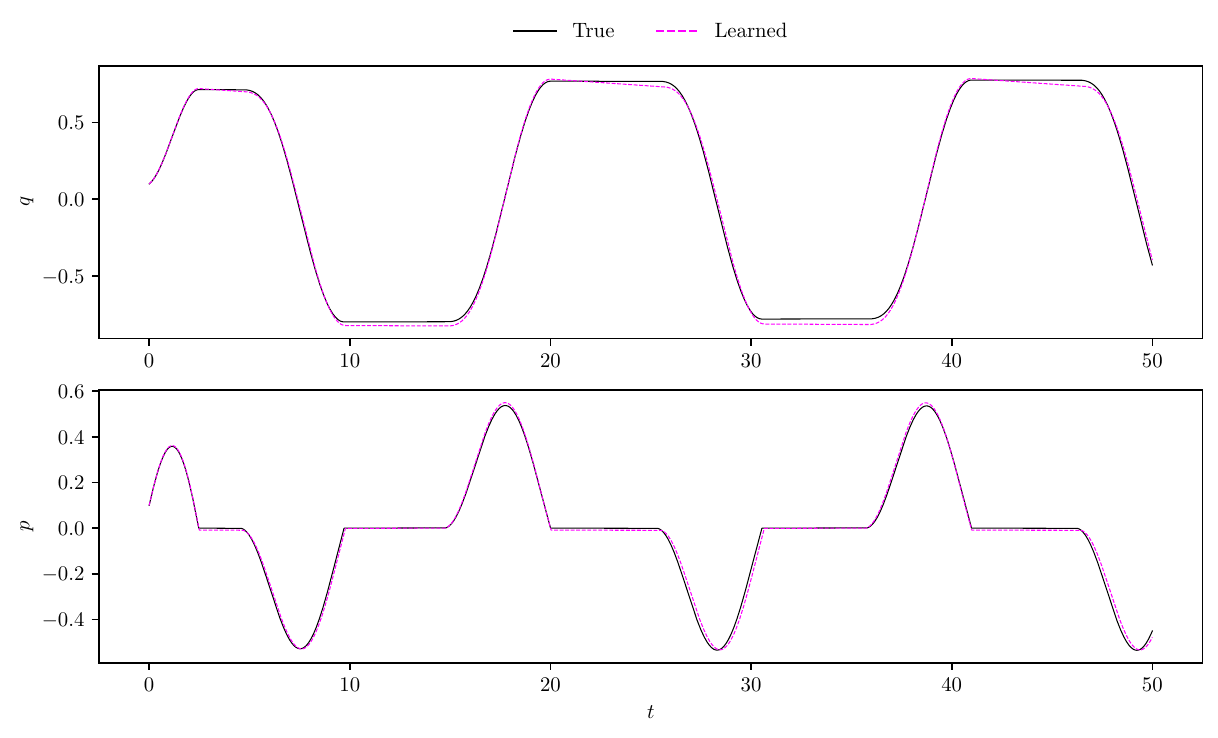}
        \caption{Case 1, recovered switch}
        \label{fig:traj-dryp-long-case1-learned-switch}
    \end{subfigure}
    \hfill
    \begin{subfigure}{0.49\linewidth}
        \centering
        \includegraphics[width=.9\linewidth]{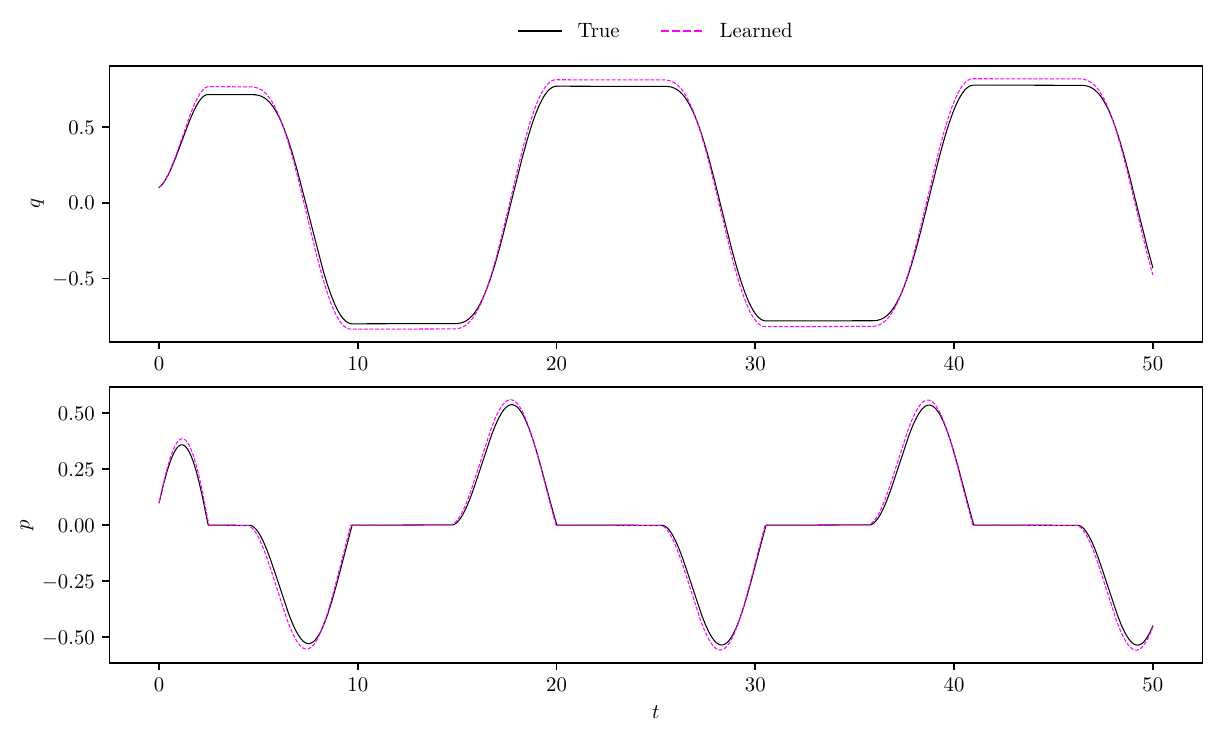}
        \caption{Case 1, known switch}
        \label{fig:traj-dryp-long-case1-known-switch}
    \end{subfigure}

    \vspace{0.75em}

    \begin{subfigure}{0.49\linewidth}
        \centering
        \includegraphics[width=.9\linewidth]{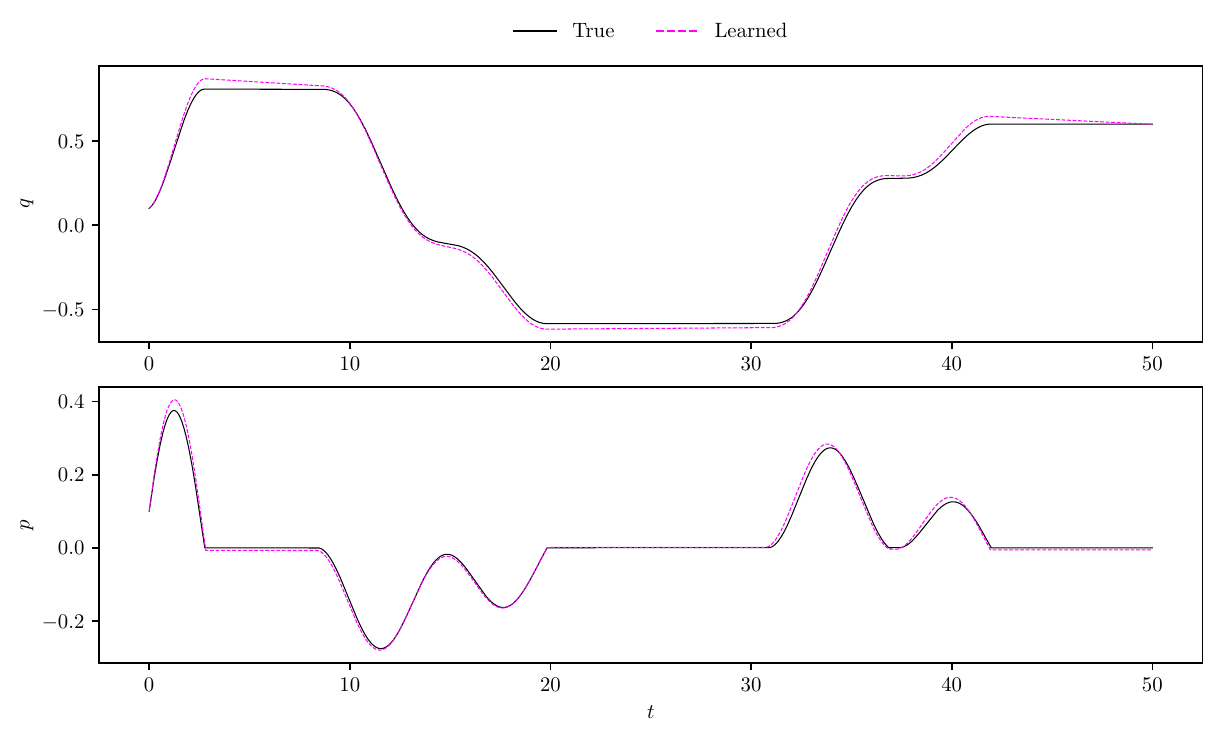}
        \caption{Case 2, recovered switch}
        \label{fig:traj-dryp-long-case2-learned-switch}
    \end{subfigure}
    \hfill
    \begin{subfigure}{0.49\linewidth}
        \centering
        \includegraphics[width=.9\linewidth]{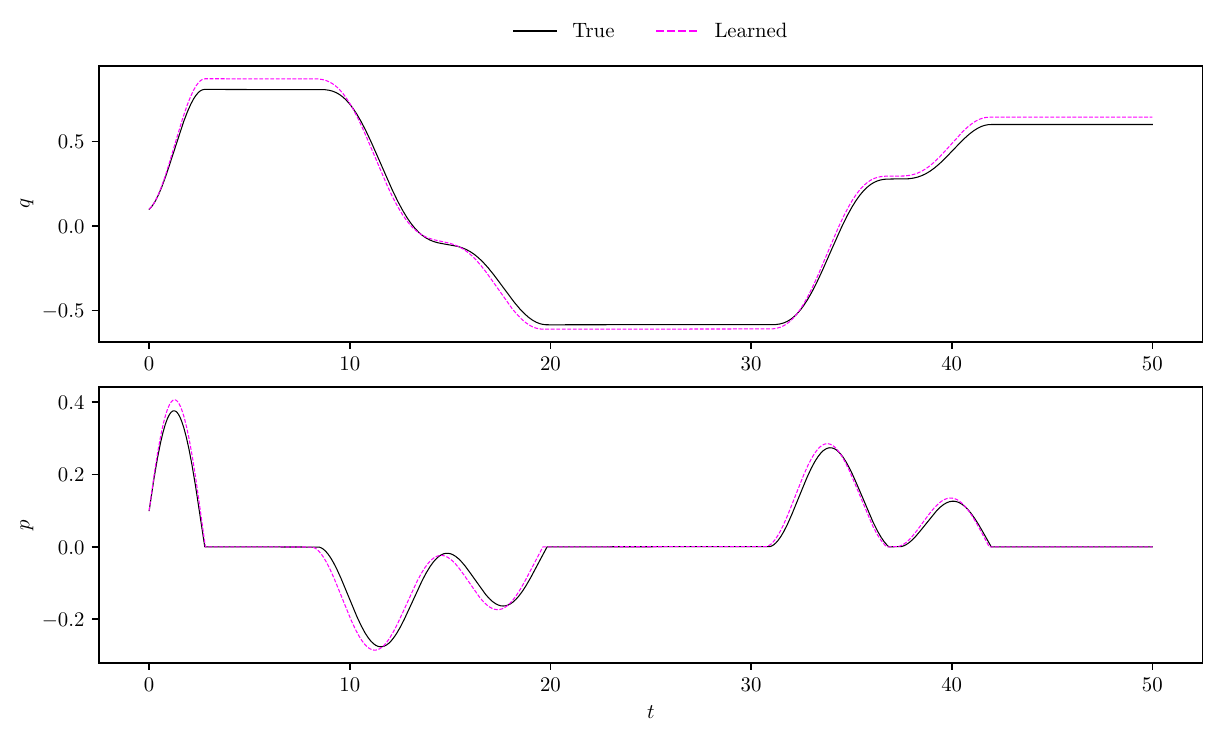}
        \caption{Case 2, known switch}
        \label{fig:traj-dryp-long-case2-known-switch}
    \end{subfigure}

    \caption{Representative long-horizon dry-friction rollouts at $\sigma=10^{-2}$.  Each panel compares the target trajectory with the learned trajectory for $q$ and $p$.  The left column shows the full
    pipeline, where the switching line is recovered from warm-up trajectories. The right column isolates the dynamics learning by supplying the true switching line.}
    \label{fig:traj-dryp-long-trajectory-known-switch-grid}
\end{figure}

The known-switch panels show that, once the switching line is supplied, the
second-order model accurately reproduces both the smooth portions and the
sliding intervals where $p=0$.  The recovered-switch panels remain close to
the target trajectories, but small discrepancies appear near transitions into
and out of sliding.  This is consistent with the main tables: the rollout
errors remain small in the mild-noise regime, while the banded vector-field errors are more sensitive to small errors in the recovered switching line.
Thus, the qualitative degradation in the full pipeline is mainly a local
switching-geometry effect rather than a global instability of the learned
dynamics.

\subsection{Unknown forcing frequency}\label{app:additionalDryp}

The dry-friction tables in the main body of the paper assume that the forcing frequency is known. To test the importance of this assumption, we repeat the strongest dry-friction configuration, namely the true switching manifold with the second-order model, but remove the known harmonic time features. The model receives random Fourier time features. The training setup, noise levels, architecture, and Filippov rollout evaluation are otherwise the same as in the main dry-friction experiments.

\begin{table}[t]
\centering
\small
\caption{Effect of removing the known forcing frequency in the dry-friction experiments. All rows use the true switching manifold and the second-order model. Known/Unknown refer to known and unknown forcing frequency $\omega$. Bold entries are the lowest between known and unknown frequency.}
\label{tab:dryp_unknown_frequency}
\begin{tabular}{lrrrrrr}
\toprule
Setting & \multicolumn{2}{c}{Rollout RMSE} & \multicolumn{2}{c}{Final RMSE} & \multicolumn{2}{c}{VF RMSE} \\
\cmidrule(lr){2-3}\cmidrule(lr){4-5}\cmidrule(lr){6-7}
 & Known & Unknown & Known & Unknown & Known & Unknown \\
\midrule
case 1, $\sigma=0$
& $\mathbf{3.18\cdot10^{-3}}$ & $2.12\cdot10^{-2}$
& $\mathbf{3.99\cdot10^{-3}}$ & $3.86\cdot10^{-2}$
& $\mathbf{6.41\cdot10^{-3}}$ & $3.47\cdot10^{-2}$ \\
case 1, $\sigma = 0.01$
& $\mathbf{2.19\cdot10^{-2}}$ & $3.47\cdot10^{-2}$
& $\mathbf{3.51\cdot10^{-2}}$ & $5.22\cdot10^{-2}$
& $\mathbf{2.91\cdot10^{-2}}$ & $5.74\cdot10^{-2}$ \\
case 1, $\sigma=0.05$
& $\mathbf{1.25\cdot10^{-1}}$ & $1.78\cdot10^{-1}$
& $\mathbf{1.84\cdot10^{-1}}$ & $2.55\cdot10^{-1}$
& $\mathbf{1.93\cdot10^{-1}}$ & $4.15\cdot10^{-1}$ \\
case 2, $\sigma=0$
& $\mathbf{2.17\cdot10^{-3}}$ & $1.69\cdot10^{-2}$
& $\mathbf{3.40\cdot10^{-3}}$ & $2.83\cdot10^{-2}$
& $\mathbf{5.89\cdot10^{-3}}$ & $4.05\cdot10^{-2}$ \\
case 2, $\sigma=0.01$
& $\mathbf{2.24\cdot10^{-2}}$ & $5.09\cdot10^{-2}$
& $\mathbf{3.56\cdot10^{-2}}$ & $8.04\cdot10^{-2}$
& $\mathbf{2.94\cdot10^{-2}}$ & $7.39\cdot10^{-2}$ \\
case 2, $\sigma=0.05$
& $\mathbf{1.26\cdot10^{-1}}$ & $2.07\cdot10^{-1}$
& $\mathbf{1.86\cdot10^{-1}}$ & $3.18\cdot10^{-1}$
& $\mathbf{1.95\cdot10^{-1}}$ & $4.47\cdot10^{-1}$ \\
\bottomrule
\end{tabular}
\end{table}

\begin{figure}[ht!]
\centering
\begin{subfigure}{0.49\linewidth}
\centering
\includegraphics[width=.9\linewidth]{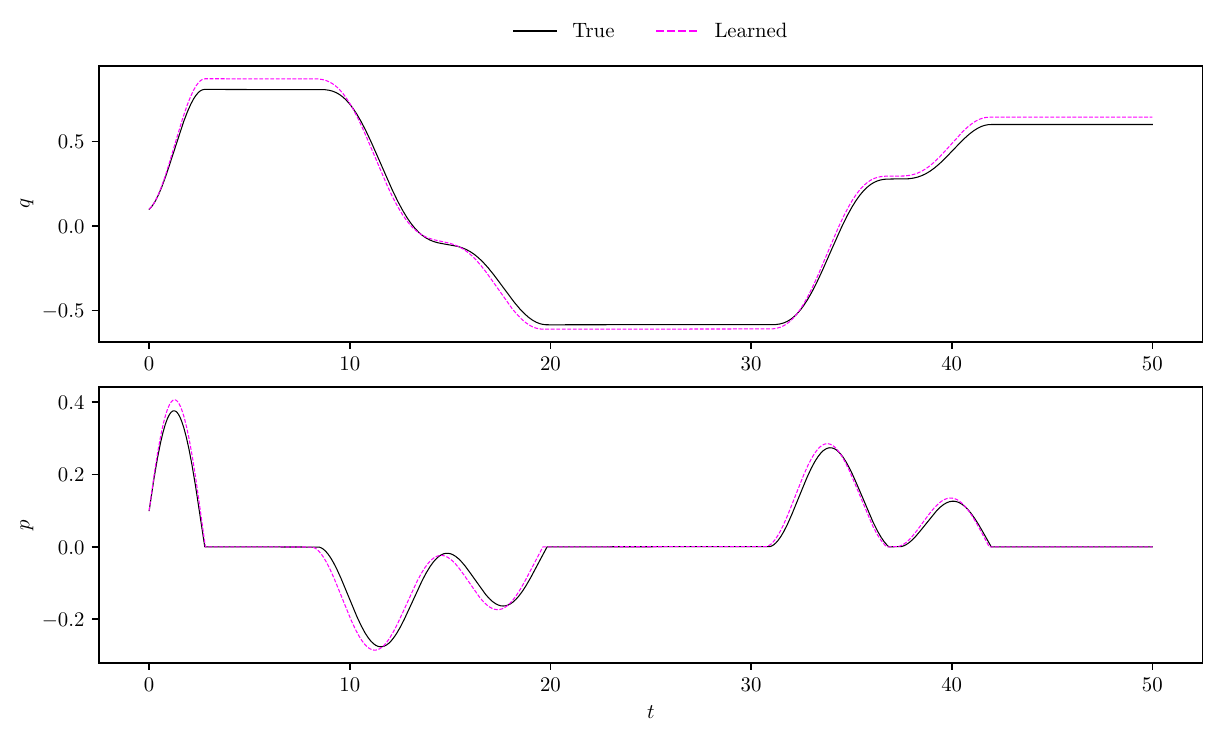}
\caption{Known forcing frequency}
\label{fig:dryp_known_frequency_appendix}
\end{subfigure}
\hfill
\begin{subfigure}{0.49\linewidth}
\centering
\includegraphics[width=.9\linewidth]{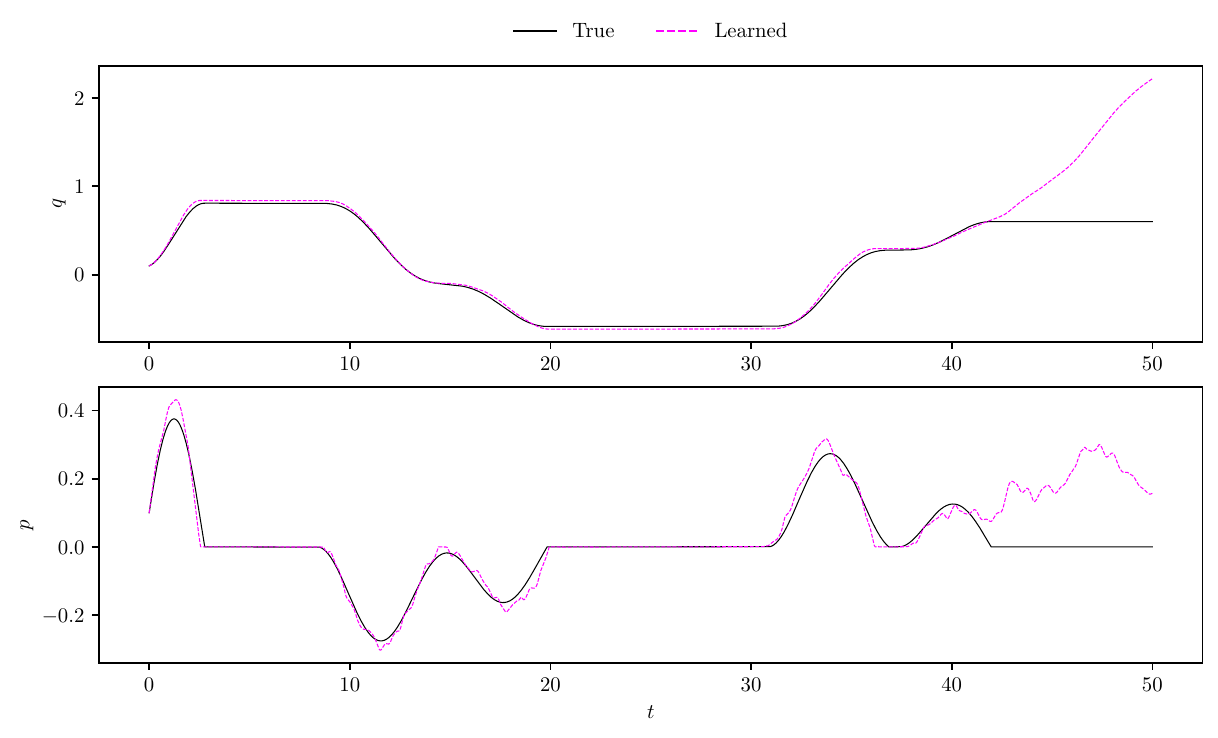}
\caption{Unknown forcing frequency}
\label{fig:dryp_unknown_frequency_appendix_panel}
\end{subfigure}
\caption{Representative long-horizon trajectories for dry-friction case 2 with \(\sigma=10^{-2}\), true switching geometry, and second-order structure. The left panel uses the known harmonic forcing frequency, while the right panel replaces it with random Fourier time features.}
\label{fig:dryp_unknown_frequency_appendix}
\end{figure}

\Cref{tab:dryp_unknown_frequency} shows that removing the known forcing frequency consistently degrades temporal extrapolation and vector-field recovery. A representative low-noise comparison is shown in \Cref{fig:dryp_unknown_frequency_appendix}. We fix the initial condition $(q_0,p_0)=(0.1,0.1)$ and $t_0=0$ for both runs. 

\section{Filippov-aware trajectory solver}
\label{app:filippovSolver}

This appendix describes the event-aware solver used for the rollout evaluations of the discontinuous models. The solver is used for single-hyperplane discontinuous models with switching function $s(x)=a^\top x+b$, and with learned one-sided vector fields \(f^+(x,t)\) and \(f^-(x,t)\).  The purpose of the solver is to avoid treating a crossing of $\Sigma:=\{x:\,s(x)=0\}$ as an
ordinary smooth step, and to evolve trajectories along \(\Sigma\) when the two
one-sided fields generate sliding.

For a step from \(x_n\) to \(x_{n+1}\) with step size \(h\), we first compute
the side \(r_n=\operatorname{sign}(s(x_n))\).  A trial step is then taken with a
branch-frozen second-order Runge--Kutta update, using \(f^+\) throughout the
step if \(r_n=+1\) and \(f^-\) throughout the step if \(r_n=-1\).  Let
\(\tilde x_{n+1}\) denote this trial endpoint.  A switching event is detected
when $s(x_n)\,s(\tilde x_{n+1})<0$. We also trigger event handling when the state is already close to the manifold,
\(|s(x_n)|<\varepsilon\), and the one-sided normal velocities indicate an
attractive sliding region.  In the experiments we use a small tolerance
\(\varepsilon\) only to decide when this special handling is needed.

When a crossing is detected, the event time is located by bisection on the
interval \([0,h]\).  During this search the same branch-frozen RK2 propagator is
used, starting from \(x_n\) on the pre-crossing side.  After bisection, the
event point is projected back onto the hyperplane,
\[
    x_\Sigma
    =
    x_{\rm event}
    -
    s(x_{\rm event})a/\|a\|_2^2,
\]
which removes the small numerical residual in the switching function.

The vector field on the switching manifold is defined using the standard
Filippov convex combination.  The trajectory is then advanced along
\(\Sigma\) using small explicit substeps and projected back onto the hyperplane
after each substep.  If the attractive condition fails, the trajectory exits the
manifold on the side indicated by the resulting one-sided normal velocities.

The event-handling procedure is repeated within one macro-step if necessary, up to a fixed small maximum number of events.  If no crossing or sliding event is detected, the solver reduces to the branch-frozen RK2 update on the current side.  This propagation rule is used for the long-horizon rollout diagnostics reported in the experiments.  The current implementation cannot handle multi-surface interactions. The event detector is also not a fully general root-finding procedure. It may miss a crossing if a trajectory touches or crosses \(\Sigma\) an even number of times within a single branch-frozen trial step, so that the endpoint sign does not change.  In the systems considered here, the chosen time steps were small enough that this event-detection rule was sufficient in practice. For more details on the implementation, see the GitHub repository \url{https://github.com/davidemurari/learn-piecewise-smooth}.

\section{Sequential recovery of multiple switching lines}
\label{app:three-line-appendix}

The geometry-recovery routine used in the main experiments is designed to find a single affine switching set from a cloud of candidate points.  This appendix records a simple extension of the same idea to a setting with more than one interface.  For this example, the data are noisy samples from three non-parallel lines in $\mathbb R^2$, plus a small number of outliers.  The purpose is to test whether the RANSAC--TLS step can be applied sequentially.

\begin{figure}[ht!]
    \centering
    \includegraphics[width=0.35\linewidth]{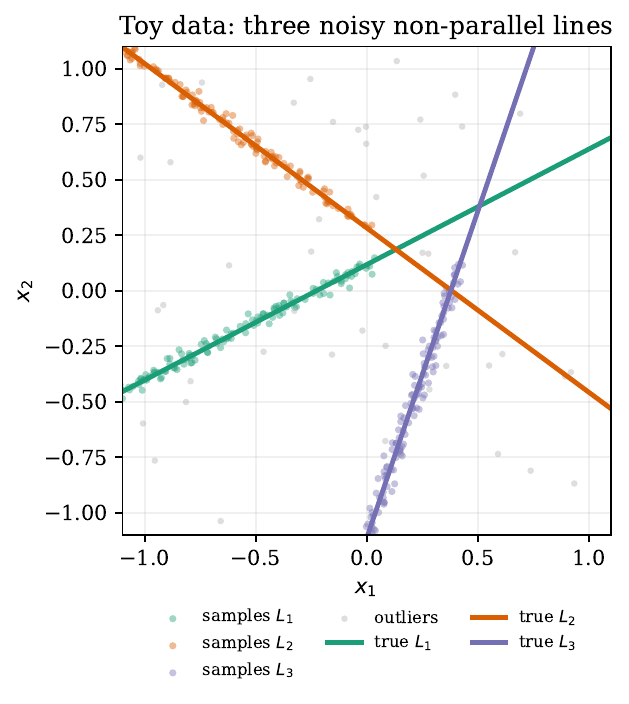}
    \caption{Toy candidate cloud used in the multi-line recovery test.  The
    coloured points are noisy samples from the three true non-parallel lines, and
    the grey points are uniform outliers.}
    \label{fig:three-line-data}
\end{figure}

\begin{figure}
    \centering
    \includegraphics[width=0.4\linewidth]{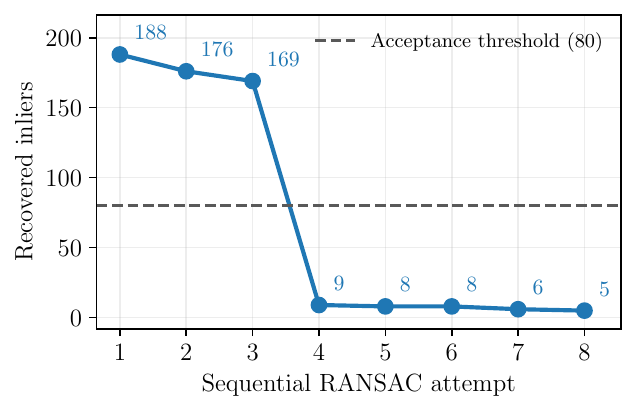}
    \caption{Number of inliers for each RANSAC attempt, after removing the previous inliers.}
    \label{fig:sequential_support}
\end{figure}

\Cref{fig:three-line-data} shows the toy data.  Each coloured cloud is sampled near one true line, with independent Gaussian perturbations.  The grey points are unrelated outliers.  For a line written as $a^\top x+b=0$,  $\|a\|_2=1$, we run the same RANSAC routine used in the geometry-learning code. After one line has been recovered, its inliers are removed from the cloud, and the procedure is repeated on the remaining points. We show in \Cref{fig:sequential_support} the number of inliers of each attempt. This figure demonstrates that when the number of hyperplanes is unknown, the number of inliers can be a viable heuristic for determining it.

\begin{figure}[ht!]
    \centering
    \includegraphics[width=.8\linewidth]{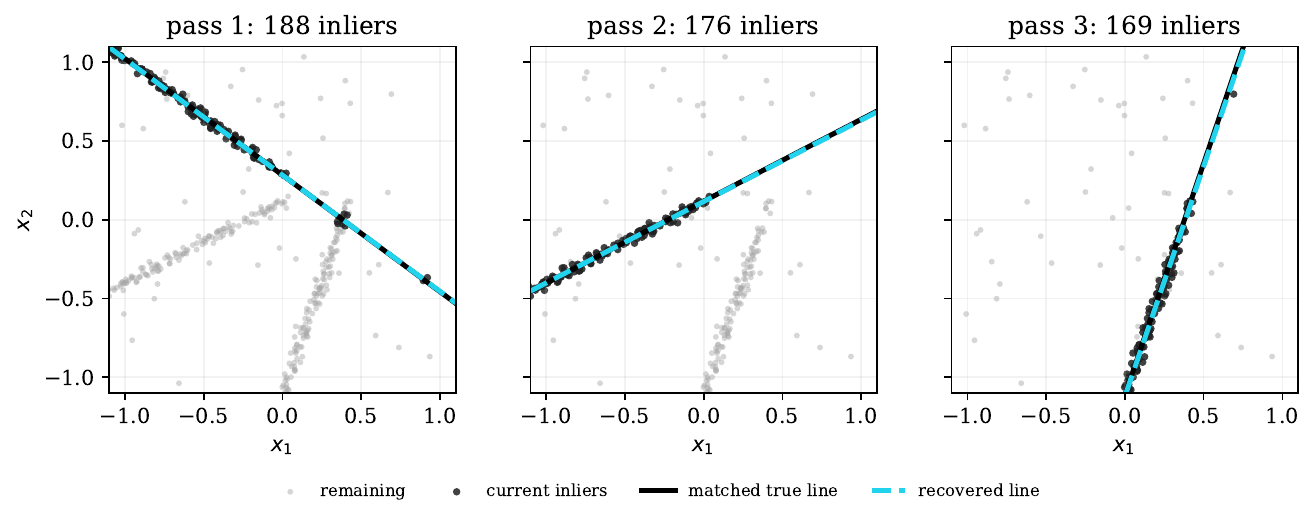}
    \caption{Sequential RANSAC--TLS extraction. In each panel, black points are
    the inliers selected during that pass, while light grey points remain available
    for later passes; points removed during earlier passes are omitted. The solid
    black line is the matched true line, and the dashed cyan line is the recovered
    TLS-refined estimate.}
    \label{fig:three-line-sequential}
\end{figure}

\Cref{fig:three-line-sequential} illustrates the three passes. The first pass matches $L_2$ with angle error $0.0854^\circ$ and offset
error $2.55\cdot10^{-4}$. The second pass matches $L_1$ with angle error
$0.135^\circ$ and offset error $1.35\cdot10^{-3}$. The third pass matches
$L_3$ with angle error $0.111^\circ$ and offset error $4.82\cdot10^{-3}$.

This toy problem also shows the main limitation of the sequential strategy.  The removal step is greedy. If the first recovered line absorbs points belonging to another nearby line, later passes inherit that mistake. Thus, this is not a complete multi-surface estimator or a full reconstruction method for polygonal cell complexes.  It is a useful diagnostic for well-separated affine interfaces, but genuinely interacting or too close switching sets would require a joint model-selection step rather than repeated single-line recovery.

\end{document}